\documentclass[11pt]{amsart}

\usepackage{latexsym,amssymb,amsmath,amscd, eucal}
\input xy
\xyoption{all}

\ifx\pdftexversion\undefined
  \usepackage[dvips]{graphics}
\else
  \usepackage[pdftex]{graphics}
\fi

\newcommand{\spec}{\mathop{\rm spec}\nolimits}

\newcommand{\im}{\mathop{\rm Im}\nolimits}

\newcommand{\Tr}{\mathop{\rm Tr}\nolimits}
\newcommand{\Ad}{\mathop{\rm Ad}\nolimits}
\newcommand{\Ker}{\mathop{\rm Ker}\nolimits}
\newcommand{\cod}{\mathop{\rm cod}\nolimits}

\newtheorem{dfn}{Definition}[section]
\newtheorem{rmk}{Remark}[section]
\newtheorem{thm}{Theorem}[section]
\newtheorem{thm*}{Theorem}
\newtheorem{cor}{Corollary}[section]
\newtheorem{prop}{Proposition}[section]
\newtheorem{lem}{Lemma}[section]

\newcommand{\Pf}{\medskip\hskip10pt\emph{Proof}. }
\newcommand{\EPf}{\hfill$\Box$\vspace{.5cm}}
\def\C{{\mathbb C}}
\def\Z{{\mathbb Z}}
\def\N{{\mathbb N}}
\def\R{{\mathbb R}}

\def\I{{\mathbf I}}
\def\L{{\mathcal L}_I}

\def\Lell{{\mathcal L}_{I,\ell}}
\def\Lellbar{\overline{\mathcal L}_{I,\ell}}

\def\O{{\mathcal O}}

\def\Vell{{\mathcal V}_{I,\ell}}

\def\W{{\mathcal W}_I}
\def\Well{{\mathcal W}_{I,\ell}}

\def\U{{\mathcal U}_I}
\def\intU{{\overset{\circ}{\mathcal U}}_I}

\def\Uell{{\mathcal U}_{I,\ell}}
\def\intUell{ {\overset{\circ}{\mathcal U}}_{I,\ell}}
\def\Uellbar{\overline{\mathcal U}_{I,\ell}}

\def\P{{\mathcal P}_I}

\def\Pell{{\mathcal P}_{I, \ell}}
\def\Pellbar{\overline{\mathcal P}_{I, \ell}}

\newcommand{\Hom}{\mathop{\rm Hom}\nolimits}
\newcommand{\LHom}{{\mathcal L}\Hom}
\newcommand{\Rep}{\mathop{\rm Rep}\nolimits}
\newcommand{\LRep}{{\mathcal L}\Rep}
\newcommand{\Sym}{{\rm Sym}_2(n)}

\newcommand{\Lg}{\mbox{$\mathfrak g$}}

\begin{document}


\title[Unitary Matrices and Lagrangian Involutions]
	{Eigenvalues of Products of Unitary Matrices and Lagrangian Involutions}

\author[Falbel]{ Elisha Falbel}

\address{Institut de Math\'ematiques \\
		Universit\'e Pierre et Marie Curie \\
		Case 82\\
		4, Place Jussieu
		F-75252 Paris, France}

\email{falbel@math.jussieu.fr}

\author[Wentworth]{ Richard A. Wentworth}

\address{Department of Mathematics \\
   Johns Hopkins University \\
   Baltimore, MD 21218}

\email{wentworth@jhu.edu}


\begin{abstract} This paper introduces a submanifold of the moduli space of unitary representations of the fundamental group of a punctured sphere with fixed local monodromy. The submanifold is defined via products of involutions through Lagrangian subspaces. We show that the moduli space of Lagrangian representations is a Lagrangian submanifold of the moduli of unitary representations.
\end{abstract}



\maketitle

\thispagestyle{empty}


\baselineskip=16pt
\setcounter{footnote}{0}


\section{Introduction}

 Let $\spec(A)$ denote the set of  eigenvalues of a unitary $n\times n$ matrix $A$.
An old problem asks the following question: 
 what are the possible collections of eigenvalues 
 $\spec(A_1),\ldots, \spec(A_\ell)$ which arise from  matrices
 satisfying $A_1\cdots A_\ell={\mathbf I}$, $\ell\geq 3\ $? 
(A review of related problems and recent developments 
can be found in  \cite{F}). 
 For an equivalent formulation in terms of representations, let
 $\Gamma_\ell$ denote the free group on $\ell-1$ generators with presentation
\begin{equation} \label{gamma}
\Gamma_\ell=\langle
 \gamma_1,\ldots, \gamma_\ell : \gamma_1\cdots\gamma_\ell =1\rangle
\end{equation}
and let $U(n)$ denote the group of 
unitary $n\times n$ matrices.  
We shall say that a collection of conjugacy classes 
$C_1,\ldots, C_\ell$ in $U(n)$  is \emph{realized} 
by a unitary representation if there is a homomorphism
 $\rho: \Gamma_\ell \to U(n)$ with $\rho(\gamma_s)\in C_s$ 
for each $s=1,\ldots, \ell$.

A natural subclass of linear representations of $\Gamma_\ell$  consists of
those generated by reflections through linear subspaces.  
In the case of unitary representations, one may consider Lagrangian planes $L$ and
 their associated involutions $\sigma_L$.   
Given a pair of Lagrangian subspaces $L_1, L_2$  in $\C^n$, the product 
$\sigma_{L_1}\sigma_{L_2}$ is an element of $ U(n)$.  
Moreover, \emph{any} unitary matrix may be obtained in this way (cf.\ Proposition
\ref{P:tau2} below).
For Lagrangians $L_1, \ldots, L_\ell$,
 one can define a unitary representation of
 $\Gamma_\ell$ via $\gamma_s\mapsto \sigma_{L_s}\sigma_{L_{s+1}}$, 
for $s=1,\ldots, \ell-1$, and $\gamma_\ell\mapsto \sigma_{L_{\ell}}\sigma_{L_1}$.
  We shall call these \emph{Lagrangian representations} (see Definition 
\ref{lhom}).  There is a natural equivalence relation obtained by rotating
every Lagrangian by an element of $U(n)$, and this corresponds to conjugation of
the representation. 
 We will say that a given  collection of conjugacy classes is 
\emph{realized} by a Lagrangian representation if the homomorphism 
$\rho$ of the previous paragraph may be chosen to be Lagrangian.

At first sight, Lagrangian representations  may seem very special.  
The main result of this paper is that in fact they exist in abundance.  We will prove

\begin{thm*}[cf.\ Section \ref{S:proof} and Propositions \ref{P:embedding}
 and \ref{maximalrank2}] \label{maintheorem}
If there exists a unitary representation  of $\Gamma_\ell$ realizing a 
given collection of conjugacy classes in  $U(n)$, then there also exists a 
Lagrangian representation realizing the same conjugacy classes. 
\end{thm*}

We also study the global structure of the moduli space of Lagrangian representations.
Let $\mathfrak a$ denote a specification 
of $\ell$ conjugacy classes $C_1,\ldots, C_\ell$, and let 
$\Rep_{\mathfrak a}^{irr.}(\Gamma_\ell, U(n))$ denote 
the set of equivalence classes of irreducible representations 
$\rho:\Gamma_\ell\to U(n)$ with each $\rho(\gamma_s)\in C_s$.  
Note that for generic choices of $\mathfrak a$, all representations are
irreducible.
Then  
$\Rep_{\mathfrak a}^{irr.}(\Gamma_\ell, U(n))$ 
is a smooth manifold which carries a symplectic structure coming from its
realization as the reduction of a quasi-Hamiltonian $G$-space (cf.\
\cite{AMM};
for a brief description, 
see Section \ref{S:symplectic}).  We refer to this as the
\emph{natural} symplectic structure.  
Let $\LRep_{\mathfrak a}^{irr.}(\Gamma_\ell, U(n))\subset  
\Rep_{\mathfrak a}^{irr.}(\Gamma_\ell, U(n))$
denote the subset of irreducible Lagrangian representations.
Then we have
\begin{thm*} \label{T:lagrangian}
With respect to the natural symplectic structure
$$\LRep_{\mathfrak a}^{irr.}(\Gamma_\ell, U(n))\subset
\Rep_{\mathfrak a}^{irr.}(\Gamma_\ell, U(n))$$
is a smoothly embedded Lagrangian submanifold.
\end{thm*}

Characterizations of which conjugacy classes are realized by products of unitary
matrices have been
 given in \cite{Be, Bi2, AW, K}.  We will give a brief review 
in Section \ref{S:inequalities} below.  The basic result is that the 
allowed region is given by a collection of affine inequalities on the log
 eigenvalues.  The ``outer walls" of the allowed region correspond to 
spectra realized only by reducible representations.  In general, there are 
also ``inner walls" corresponding to spectra that are realized by both
 reducible and irreducible representations.  The open chambers
 complementary to these walls correspond to spectra that are
 realized only by irreducible representations.  The term ``generic"
 used above
 refers to spectra in the open chambers.

This structure suggests a proof of 
Theorem \ref{maintheorem} via induction on the rank and
 deformation theory, and this is the approach we shall take.  
In Section \ref{S:lagrangian}, we prove some elementary facts about configurations
of pairs and triples of Lagrangian subspaces in $\C^n$.  We define Lagrangian
representations and 
discuss their relationship to unitary representations. 
In particular, we show that the Lagrangian representation
space is isotropic with respect to the natural symplectic
structure.
In Section \ref{S:deformation}, 
 after briefly
reviewing the case of unitary representations, 
we develop the deformation theory of Lagrangian
representations in more detail.
We introduce two methods to produce a family of Lagrangian representations from a
given one.  We call these deformations
 \emph{twisting} and \emph{bending} (see Definitions
\ref{D:twist} and \ref{D:bending}), and they are in part
motivated by the geometric flows studied by Kapovich and Millson \cite{KM}. 
We prove that twisting and bending deformations, applied to an irreducible
Lagrangian representation, span all possible variations of the conjugacy classes (see
Proposition \ref{maximalrank2}).  
As a consequence, if there is a single point interior to one of the  
chambers described above that is  realized by a Lagrangian representation, 
then all points in the chamber are also realized by Lagrangians  
 (see Corollary \ref{isolatedchambers}).   
This reduces the existence problem to ruling out the possibility of isolated chambers 
realized by unitary representations, but not by Lagrangians.  
To achieve this we make a detailed analysis of the wall structure in Section
\ref{S:proof}.  
A basic fact is that any reducible Lagrangian representation may be 
perturbed to an irreducible one.  Hence, inductively, any chamber having an 
outer wall as a face is necessarily populated by Lagrangian representations.
A topological argument that exploits an estimate (Proposition
\ref{P:codimension})  on the codimension
of the set of reducible representations 
shows that inner walls may also be ``crossed" by Lagrangian representations.

It should be apparent from this description that 
 our proof of Theorem \ref{maintheorem} is somewhat indirect. 
 A more precise description of the obstructions to deformations of 
reducible unitary and Lagrangian representations is desirable. 
In \cite{FMS} Lagrangians were used to give a geometrical explanation
 of the inequalities for $U(2)$ representations in terms of 
spherical polygons.  For higher rank  it is tempting to look for a
 similar geometrical interpretation of the inequalities,
 though we have not obtained such at present.
Unitary representations of surface groups are related to stability
 of holomorphic  vector bundles through the famous theorem of 
Narasimhan and Seshadri \cite{NS} and its generalization to punctured
 surfaces by Mehta and Seshadri \cite{MS}.  A challenging problem is to
 give an analytic description of those holomorphic structures
 which give rise to Lagrangian representations.

 We conclude this introduction by pointing out
 an alternative interpretation of the result in Theorem \ref{maintheorem}.
 Let us say that matrices $A_1,\ldots, A_\ell\in U(n)$ are 
\emph{pairwise  symmetrizable} if for each $s=1,\ldots,\ell$, there is 
$g_s\in U(n)$ so that both $g_s A_s g_s^{-1}$ and $g_s A_{s+1}g_s^{-1}$ are 
symmetric  (where $A_{\ell+1}=A_1$).  Also, throughout the paper, for 
unitary matrices $A$ and $B$, $A\sim B$ indicates that $A$ and $B$ are
conjugate.  We then have the following reformulation of Theorem \ref{maintheorem}.
\begin{thm*} \label{symmetrizable}
  Given $n\times n$ unitary matrices $\{A_s\}_{s=1}^\ell$, $A_1\cdots A_\ell=\I$, there exists a possibly different collection of unitary matrices $\{B_s\}_{s=1}^\ell$, $B_1\cdots B_\ell=\I$, 
$A_s\sim B_s$ for $s=1,\ldots, \ell$, such that $B_1,\ldots, B_\ell$ are pairwise symmetrizable.
\end{thm*}
See
Section \ref{S:lrep} for the proof.


\section{Unitary Representations} \label{S:unitary}

\subsection{The space of conjugacy classes} \label{S:conjugacy}

  We begin with some notation. 
Given integers $n\geq 1$ and  $\ell\geq 3$:
\begin{itemize}
\item Let $\overline{\mathcal M}_\ell(n)$ denote the set of all
  $\ell\times n$ matrices ${\mathfrak a}=(\alpha^s_j)$,
 $1\leq s\leq\ell$, $1\leq j\leq n$,
 where for each $s$,  $\alpha^s=(\alpha^s_1, \ldots,\alpha^s_n)$
 satisfies $0\leq \alpha^s_1\leq \cdots \leq \alpha^s_n\leq 1$.  
 \item Let $\overline{\mathcal A}_\ell(n)$ be the quotient of
 $\overline{\mathcal M}_\ell(n)$ defined by the following equivalence:
  identify a point of the form $\alpha^s=
 (\alpha^s_1, \ldots,\alpha^s_k, 1, \ldots, 1)$, $\alpha^s_k<1$,
  with $\tilde\alpha^s=
 (0,\ldots, 0, \tilde\alpha^s_{n-k+1}, 
\ldots,\tilde\alpha^s_n)$, where $\tilde\alpha^s_{n-k+i}=\alpha^s_i$, 
$i=1,\ldots, k$. 
 \item
Let ${\mathcal A}_\ell(n)\subset \overline{\mathcal A}_\ell(n)$ be the
 open subset where all inequalities are strict: 
$0< \alpha^s_1< \cdots < \alpha^s_n< 1$, for each $s$.
\end{itemize}
   
 \noindent  For each ${\mathfrak a}\in\overline{\mathcal A}_\ell(n)$
 we define the \emph{index} as follows: 
 choose the representative of $\mathfrak a$ where
 $0\leq \alpha^s_1\leq \cdots \leq \alpha^s_n< 1$, for each $s$,  and set
\begin{equation} \label{index}
I({\mathfrak a})=\sum_{s=1}^\ell \sum_{j=1}^n\alpha^s_j\ .
\end{equation}
We define 
$\overline{\mathcal A}_\ell^{\Z}(n)=\{ {\mathfrak a}\in 
\overline{\mathcal A}_\ell(n) : I({\mathfrak a})\ \text{is an integer}
\, \}$, 
${\mathcal A}_\ell^{\Z}(n)={\mathcal A}_\ell(n)\cap
\overline{\mathcal A}_\ell^{\Z}(n)$.

\begin{dfn} \label{mplane}
  For a nonnegative integer $I$, define the \emph{open M-plane} 
  by  
$$
\Pell(n)=\{ {\mathfrak a}\in {\mathcal A}_\ell^{\Z}(n) : I({\mathfrak a})=I \}\ .
  $$
   The closure $\Pellbar(n)$ of $\Pell$ in 
$\overline{\mathcal A}_\ell^{\Z}(n)$ will be called the 
\emph{closed M-plane}.  Finally, let
$$\Pellbar^\ast(n)=\{ {\mathfrak a}\in \Pellbar : I({\mathfrak a})=I \}\ .$$
\end{dfn}
 
 \noindent  
Observe that $\Pellbar(n)$ is  a closed connected cell.
 Notice also that the closed $M$-planes are not disjoint,
 whereas of course $\Pellbar^\ast(n)\cap 
\overline{\mathcal P}_{J,\ell}^\ast(n)=\emptyset $ if $I\neq J$.
  We therefore have a \emph{disjoint} union
$$
 \overline{\mathcal A}_\ell^{\Z}(n)=
\bigcup_{0\leq I \leq n\ell-1} \Pellbar^\ast(n)\ .
$$

 For each $s$ choose a partition $m^s$ of 
$\{1,\ldots, n\}$, i.e.\ a set of integers 
$0=m^s_0<m^s_1<\cdots<m^s_{l_s}=n$. 
Here, $l_s$ is the \emph{length} of the partition. 
Specifying $l_s$ numbers $0\leq\hat\alpha^s_1<\ldots 
<\hat\alpha^s_{l_s}<1$ along with a partition of length 
$l_s$ uniquely determines a point in 
${\mathfrak a}=(\alpha^s_j)\in\overline{\mathcal A}_\ell(n)$,
 where 
 $\alpha^s_i=\hat\alpha^s_j$ for $m^s_{j-1}<i\leq m^s_{j}$.
 Conversely, given a point
$ {\mathfrak a}\in\overline{\mathcal A}_\ell(n)$
 with the distinct entries $0\leq\hat\alpha^s_1<\ldots
 <\hat\alpha^s_{l_s}<1$, a partition of length $l_s$ 
is determined by the multiplicities 
$\mu_j^s=m^s_j-m^s_{j-1}$  of the $\hat\alpha^s_j$. 
We shall say that $\alpha^s$ \emph{has  the multiplicity 
structure of} $m^s$.  

Let ${\mathfrak m}=(m^1,\ldots, m^\ell)$ be a choice of $\ell$
 partitions.  In addition, choose a (possibly empty subset)
 $z\subset\{1,\ldots, \ell\}$ of cardinality $|z|$.
This data leads to  the following refinement of the $M$-plane.
\begin{align*} 
\Pell(n, {\mathfrak m}, z)&=\bigl\{ {\mathfrak a}=(\alpha^s_j)\in 
\Pellbar^\ast(n) : \ \alpha^s \ \text{has multiplicity structure}\ 
m^s\ \text{for all}\ s\ , \\
&\qquad\qquad\qquad \text{and}\ \hat\alpha_1^s=0\
 \text{if and only if}\ s\in z\ \bigr\}\ ; \\
\Pellbar(n, {\mathfrak m}, z)&=
\ \text{the closure of}\ \Pell(n, {\mathfrak m}, z) \ \text{in }\ 
\overline{\mathcal A}_\ell^{\Z}(n) \ ; \\
\Pellbar^\ast(n, {\mathfrak m}, z)&=\Pellbar(n, 
{\mathfrak m}, z)\cap\Pellbar^\ast(n)\ . \notag
\end{align*}

 Next,  notice that there is a natural partial ordering on
 multiplicities:
if ${\mathfrak p}=(p^1,\ldots, p^\ell)$ and ${\mathfrak m}=(m^1,\ldots, m^\ell)$, we say that
 ${\mathfrak p}\leq {\mathfrak m}$ if for each $s=1,\ldots, \ell$ the partition $p^s$ is a subset of $m^s$.
We then have a stratification by the cells $\Pell(n,{\mathfrak m},z)$  in the sense that
$$
 \Pellbar^\ast(n,{\mathfrak m},z) =\bigcup_{{\mathfrak p}\leq 
{\mathfrak m}\, ,\, z\subset \tilde z\subset\{1,\ldots, \ell\}} 
\Pell(n,{\mathfrak p},\tilde z)\ .
$$
  In particular,
$$
 \Pellbar^\ast(n) =\bigcup_{  
\mathfrak m \, ,\, z\subset\{1,\ldots,\ell\}     }
    \Pell(n,{\mathfrak m},z)
$$

There is a similar, though slightly more complicated, stratification of $\Pellbar(n,{\mathfrak m},z)$
which involves strata of lower index. 
 To describe this, consider the limit $\bar{\mathfrak a}$ in
 $\overline{\mathcal A}_\ell^{\Z}(n)$ of points in
 $\Pell(n,{\mathfrak m},z) $ where $\hat\alpha^{s_0}_{l_{s_0}}\to 1$, 
for some $s_0\in\{1,\ldots, \ell\}$, but the $\hat\alpha^s_{l_s}$
 remain bounded away from $1$ for $s\neq s_0$.  From the
  defining equivalence $\overline{\mathcal M}_\ell(n)\to
 \overline{\mathcal A}_\ell(n)$ and the convention \eqref{index}
 for the index, it follows that
$$
\overline I=I(\bar{\mathfrak a})= I-(n-m^{s_0}_{l_{s_0}-1})<I\ .
$$
Furthermore, we may define a new collection of partitions $\bar{\mathfrak m}$, $\bar m^s(\bar l_s)=m^s(l_s)$ for $s\neq s_0$, and
\begin{align*}
\text{ if}\ s_0\in z\, ,\ \text{then}\
&\begin{cases}
\bar m^{s_0}_i = m_i^{s_0} + (n-m^{s_0}_{l_{s_0}-1})\ ,\ 1\leq i\leq l_{s_0}-1\ ,\\
\bar l_{s_0} = l_{s_0} -1\ , \\
\bar z=z\ ;
\end{cases} \\
\text{if}\ s_0\not\in z\, ,\ \text{then}\
&\begin{cases}
\bar m^{s_0}_1 =n-m^{s_0}_{l_{s_0}-1}\ ,\\
\bar m^{s_0}_{i+1} = m_i^{s_0} + (n-m^{s_0}_{l_{s_0}-1})\ ,\ 1
\leq i\leq l_{s_0}-1\ ,\\
\bar l_{s_0} = l_{s_0} \ , \\
\bar z=z\cup\{s_0\}\ .
\end{cases} 
\end{align*}
With these definitions, it is clear that $\bar{\mathfrak a}
\in{\mathcal P}_{\overline I, \ell}(n,\bar{\mathfrak m},\bar z)$.
  A  stratification of $\Pellbar(n,{\mathfrak m},z)$ is then 
obtained by adding, in addition to sets of the form 
$\Pell(n, {\mathfrak p}, \tilde z)$, all sets
 ${\mathcal P}_{\overline I, \ell}(n,\bar{\mathfrak m},\bar z)$
 derived from these strata in the manner described above.


\subsection{Inequalities for unitary representations}
 \label{S:inequalities}

Let $\Gamma_\ell$ be as in \eqref{gamma}, and fix an integer $n\geq 1$.
 We will denote the $U(n)$-\emph{representation variety} of 
$\Gamma_\ell$ by
$$
\Hom(\Gamma_\ell, U(n))=\{\text{homomorphisms}\ \rho: \Gamma_\ell \to U(n) \}\ .
$$
  We denote the subspaces of irreducible and reducible homomorphisms by $\Hom^{irr.}(\Gamma_\ell, U(n))$ and $\Hom^{red.}(\Gamma_\ell, U(n))$, respectively.  
 The group $U(n)$ acts on $\Hom(\Gamma_\ell, U(n))$ (say, on the left) 
by conjugation.  We define the \emph{moduli space of representations} to be the
quotient
$$
\Rep(\Gamma_\ell, U(n))=U(n)\bigr\backslash \Hom(\Gamma_\ell, U(n))\  .
$$
Following the notation for homomorphisms,  
subsets of equivalence classes of irreducible and reducible 
homomorphisms are denoted  by $\Rep^{irr.}(\Gamma_\ell, U(n))$
 and $\Rep^{red.}(\Gamma_\ell, U(n))$, respectively. 
  With the presentation of $\Gamma_\ell$ given in \eqref{gamma}, to
 each $[\rho]\in \Rep(\Gamma_\ell, U(n))$ we associate 
 conjugacy classes $\rho(\gamma_1), \ldots, \rho(\gamma_\ell)$.  
In this section, we give a brief description of which collections of
  $\ell$ conjugacy classes are realized by unitary representations in
 this way.
 
Given $A\in U(n)$, we may express its eigenvalues as 
$(\exp(2\pi i\alpha_1), \ldots, \exp(2\pi i\alpha_n))$,
 with $0\leq\alpha_1\leq \cdots\leq \alpha_n< 1$, and this
 expression is unique.  We will therefore write: 
$\spec(A)=\alpha=(\alpha_1,\ldots, \alpha_n)$. 
 The spectrum determines and is determined uniquely by the 
conjugacy class of $A$. 
If $A_1,\ldots, A_\ell\in U(n)$, 
$A_1\cdots A_\ell={\mathbf I}$, and $\spec(A_s)=\alpha^s$, 
then by taking determinants
 we see that the index $I(\alpha^s_j)$ defined in \eqref{index}
is an \emph{integer}.  
As in the introduction, we may recast this in terms of representations.
For $\rho\in \Hom(\Gamma_\ell, U(n))$, we  set 
$A_s=\rho(\gamma_s)$, and 
 there is a well-defined integer $I=I(\rho)$
 associated to $\rho$.  Clearly, $I(\rho)$ depends only on the
 conjugacy class of the representation, so it is actually well-defined
 for $[\rho]\in \Rep(\Gamma_\ell, U(n))$.  

\begin{dfn} \label{specproj}
Given $\rho\in \Hom(\Gamma_\ell, U(n))$, the integer
 $I(\rho)$ is called the \emph{index} of the representation.  
We define the \emph{spectral projection}
$$
\pi: \Hom(\Gamma_\ell, U(n))\longrightarrow 
\overline{\mathcal A}_\ell^{\Z}(n)\ :
 \quad
\rho\longmapsto [\spec(\rho(\gamma_1)),\ldots, \spec(\rho(\gamma_\ell))]\ .
$$
Then $\pi$ factors through a map {\rm (}also denoted $\pi${\rm )}
 on $\Rep(\Gamma_\ell, U(n))$. 
We denote the fibers of $\pi$ over 
${\mathfrak a}\in \overline{\mathcal A}_\ell^{\Z}(n)$ by
\begin{align*} 
\Hom_{\mathfrak a}(\Gamma_\ell, U(n))&=\pi^{-1}({\mathfrak a})\subset
 \Hom(\Gamma_\ell, U(n)) \\
\Rep_{\mathfrak a}
(\Gamma_\ell, U(n))&=\pi^{-1}({\mathfrak a})\subset \Rep(\Gamma_\ell, U(n))  
 \ . 
\end{align*}
\end{dfn}
  The image of $\pi$ is our main focus in this section.

\begin{dfn}  \label{D:unitary}
Let $\Uellbar^\ast(n)= \pi(\Hom(\Gamma_\ell, U(n)))\cap \Pellbar^\ast(n)$. 
For each collection of multiplicities
${\mathfrak m}=(m^s)$ and subsets $z\subset\{1,\ldots, \ell\}$,
we set
$$\Uell(n,{\mathfrak m}, z)=\Uellbar^\ast(n)\cap\Pell(n,{\mathfrak m}, z)\ .$$
\end{dfn}

\begin{dfn} \label{nondegenerate}
Denote the interior points of $\Uell(n,{\mathfrak m},z)$ in $\Pell(n,{\mathfrak m},z)$
 by $\intUell(n, {\mathfrak m},z)$.  A stratum $\Pell(n,{\mathfrak m},z)$ is called
\emph{nondegenerate} if either
$$\Uell(n,{\mathfrak m},z)=\emptyset\ ,$$
  or
$$\intUell(n, {\mathfrak m},z)\neq\emptyset\ .$$
\end{dfn}

The regions $\Uell(n, {\mathfrak m},z)$ have the following simple 
description (cf.\ \cite[Theorem 3.2]{Bi2} and \cite{Be, AW, K}).
\begin{thm} \label{structure}
There is a finite 
collection $\Phi_{I,\ell}(n)$ of affine linear 
functions of the $\{\alpha^s_j\}$  such that
$$
\Uellbar^\ast(n)=\left\{ {\mathfrak a}\in
\Pellbar^\ast(n) : \phi({\mathfrak a})\leq 0\ 
\text{ for all } \phi\in\Phi_{I,\ell}(n)\ \right\}\ .
$$
Moreover, the sets $\Phi_{I,\ell}(n)$, as $I$ varies,  are 
compatible with the stratification described in the previous section.
\end{thm}

\begin{dfn}  \label{outerwall}
For each $\phi\in\Phi_{I,\ell}(n)$ we define the
 \emph{outer wall} associated to $\phi$ by
$$
 W_\phi=\left\{ {\mathfrak a}\in\Pell(n,{\mathfrak m},z) : 
\phi({\mathfrak a})= 0\ \right\}\ .
 $$
 We denote the union of all outer walls by
 $$
 \Well(n,{\mathfrak m},z)=\bigcup_{\phi\in\Phi_{I,\ell}(n)} W_\phi\ .
 $$
 \end{dfn}
 
 \noindent   It follows that $\Uell(n,{\mathfrak m},z)$ is the 
closure in $\Pell(n,{\mathfrak m},z)$ of a convex connected component of 
 $\Pell(n,{\mathfrak m},z)\setminus \Well(n,{\mathfrak m},z)$.  
The representations with $\pi(\rho)\in\Well(n,{\mathfrak m},z)$ are
 reducible (see Proposition \ref{P:structure}). 
 Indeed, the functions $\phi$ defining the walls
are all of the following type.
   Fix an integer $1\leq k<n$.  Choose $\wp_{(k)}=(\wp_{(k)}^1,\ldots,
\wp_{(k)}^\ell)$, where
 for each $s=1,\ldots, \ell$,  $\wp_{(k)}^s$ is a subset of $\{1,\ldots, n\}$ of cardinality $k$.  We define a relative index by 
 \begin{equation} \label{relativeindex}
 I({\mathfrak a}, \wp_{(k)})=\sum_{s=1}^\ell\sum_{\alpha^s_j\in\wp^s_{(k)}}\alpha^s_j\ .
 \end{equation}
 Notice that  for ${\mathfrak a}\in\Uellbar^\ast(n)$
 the value of $ I({\mathfrak a}, \wp_{(k)})$ may 
a priori be any real number less than $I$.
 Suppose $\rho\in \Hom_I(\Gamma_\ell,U(n))$ is reducible.  Hence, there is a reduction $\rho: \Gamma_\ell\to U(k)\times U(n-k)$ for some $1\leq k<n$.  
The set of eigenvectors of $\rho(\gamma_s)$ lying in the $U(k)$ factor
 gives a collection of subsets $\wp_{(k)}^s$.  
Moreover, it follows, again by taking determinants that the relative 
index $I(\pi(\rho), \wp_{(k)})$ is equal to some integer $K$,
 $0\leq K\leq I$.   We will say that the reducible representation is
 \emph{compatible} with $(K,\wp_{(k)})$ if the pair
$(K, \wp_{(k)})$
  arises from some reduction of $\rho$.
The functions $\phi\in\Phi_{I,\ell}(n)$ are all of the form 
 $\phi({\mathfrak a})=I({\mathfrak a}, \wp_{(k)})-K$, for various
 choices of partitions $\wp_{(k)}$ and integers $K$.

 It is not necessarily the case, however,  that every reducible
 $\rho$ projects via $\pi$  to an \emph{outer} wall. 
Nevertheless, we see that  
there is still a hyperplane associated to any reducible. 
 This motivates the following
 
 \begin{dfn}  \label{innerwall}
 Let $\Psi_{I,\ell}(n)$ be the finite collection of affine linear functions of the form $\psi({\mathfrak a})=
  I({\mathfrak a}, \wp_{(k)})-K$,  for partitions $\wp_{(k)}$ and positive integers $K$, such that
there is some reducible $\rho$ compatible with $(K, \wp_{(k)})$ 
for which $\pi(\rho)\in\intUell(n,{\mathfrak m},z)$, for some $\mathfrak m$, $z$.
 For $\psi\in \Psi_{I,\ell}(n)$  we define the \emph{inner wall} associated to 
$\psi$ by
 $$
 V_\psi=\left\{ {\mathfrak a}\in\Pell(n,{\mathfrak m},z) : \psi({\mathfrak a})= 0\ \right\}\ .
 $$
 We denote the union of all inner walls by
 $$
 \Vell(n,{\mathfrak m},z)=\bigcup_{\psi\in\Psi_{I,\ell}(n)} V_\psi\ .
 $$
 \end{dfn}

\noindent   Hence, the distinction between the two types of walls 
is that there are points of $\Uell(n,{\mathfrak m},z)$ 
on either side of an inner wall,
 whereas $\Uell(n,{\mathfrak m},z)$ lies on only one side of each outer wall.

The precise determination of the functions in 
$\Phi_{I,\ell}(n)$  is quite involved.  
In Section \ref{S:examples}, we give the result for $\Phi_{I,3}(2)$ and $\Phi_{I,3}(3)$.
One way to view the origin of  these conditions is via the notion of
 stable and semistable parabolic structures on holomorphic vector 
bundles over $\C P^{1}$.  We will require very few details of 
this theory; the interested reader may consult  the references cited 
above. 
The following two results are consequences of this 
holomorphic description. First, we have

\begin{prop} \label{P:structure}
Let  $\rho\in \Hom_I(\Gamma_\ell,U(n))$ with $\pi(\rho)\in 
\Pell(n,{\mathfrak m},z)$.
\begin{enumerate}
\item If $\pi(\rho)\in\Well(n,{\mathfrak m},z)$, then $\rho$ is
 reducible.
\item If $\rho$ is reducible, then $\pi(\rho)\in 
\Well(n,{\mathfrak m},z)\cup\Vell(n,{\mathfrak m},z)$.
\item 
 If $\pi(\rho)\in\intUell(n, {\mathfrak m},z)$, 
there is an irreducible  representation
$ \tilde\rho$ with $\pi( \tilde\rho)={\mathfrak a}$.
\end{enumerate}
\end{prop}

\Pf Part (1) follows from the fact that an irreducible representation corresponds
to a stable parabolic structure. And if a parabolic structure is stable for a
given set of weights, it is also stable for a sufficiently small  neighborhood of
weights (an alternative, purely representation theoretic
  proof of this follows from the arguments in Section
\ref{S:deformation} below). Part (2) is by definition.  
Part (3) is immediate from \cite[Theorem 3.23]{Bi2}, since if the strict
inequalities are satisfied there exists a stable parabolic structure.
Stable structures, as mentioned, correspond to irreducible representations.
\EPf

Next, we give sharp bounds on the index.

\begin{thm} \label{indexbounds1}
For any representation $\rho: \Gamma_\ell\to U(n)$ we have
$$
n-N_0(\rho)\leq I(\rho)\leq n(\ell-1)+N_0(\rho)-N_1(\rho)\ ,
$$
where $N_0(\rho)$ is the number of trivial representations appearing in the
decomposition of $\rho$ into irreducibles, and $N_1(\rho)$ is the total
multiplicity of the eigenvalue $0$ among $\alpha^s=\rho(\gamma_s)$ for all  $s=1, \ldots,
\ell$.  Moreover, these bounds are sharp.
\end{thm}

\Pf   
The case $n=1$ is straightforward.  For $n\geq 2$, we first show
that $I(\rho)\geq n-N_0(\rho)$.  Since both sides of this inequality are
additive on reducibles, an inequality
$I(\rho)\geq n$ for irreducible representations proves the result in general by
induction.  Hence, suppose
$\rho: \Gamma_\ell\to U(n)$ is an irreducible representation with
$\pi(\rho)=(\alpha^s_j)$ and $I(\rho)<n$.
Associated to $\rho$ is a stable parabolic bundle on $\C P^{1}$ with weights
$(\hat\alpha^s_j)$ whose underlying holomorphic bundle $E$ has degree $
-I(\rho)$ (cf.\ \cite{MS}).  By the well-known theorem of Grothendieck, $E\to
\C P^{1}$ is holomorphically split into a sum of line bundles: $E=\O(d_1)\oplus
\cdots \oplus \O(d_n)$, where $\O(d)$ denotes the (unique up to isomorphism)
holomorphic line bundle of degree $d$ on $\C P^{1}$.  By assumption
$ \sum_{j=1}^n d_j=\deg E=-I(\rho)>-n $.
Hence, there is some $d_j\geq 0$.  But then $E$ contains a subbundle $\O(d_j)$
with nonnegative parabolic degree.  This contradicts  parabolic stability,
and hence also the assumption $I(\rho)<n$.  
Thus, the inequality $I(\rho)\geq n$ for irreducibles holds.
Next, notice that to any representation $\rho:\Gamma_\ell\to U(n)$ we
may associate a dual representation $\rho^\ast:\Gamma_\ell\to U(n)$
defined by: $\rho^\ast(\gamma_s)=\rho(\gamma_{\ell+1-s})^{-1}$, $s=1, \ldots,
\ell$. Using the convention \eqref{index} it follows that
$ I(\rho^\ast)=n\ell -I(\rho)-N_1(\rho) $,
where $N_1(\rho)$ is defined in the statement of the theorem.  Combining
this with the previous result $I(\rho)\geq n$, we see that
$I(\rho)\leq n(\ell-1)-N_1(\rho)$, for $\rho$ irreducible.  This argument generalizes to
the case where $\rho$ contains  trivial factors as well.  This completes the
proof of the inequality.  To prove
that the bounds are sharp we need only remark that both sides of the
inequalities are additive on reducibles and that the bounds are
evidently sharp for the case $n=1$. 
\EPf

In Section 3, we will indicate a ``Lagrangian" proof of this result
for the case $\ell=3$ (see Proposition \ref{indexbounds2}).
We conclude this section with one more

\begin{dfn} \label{chambers}
A connected component  of 
$$\Uell(n,{\mathfrak m},z)\setminus \left\{\Well(n,{\mathfrak m},z)\cup\Vell(n,{\mathfrak m},z)\right\}$$
 will be called a \emph{chamber}.  
\end{dfn}

\begin{rmk} \label{chamberstructure}
\begin{enumerate}
\item From the description given above the chambers 
 of $\Pell(n,{\mathfrak m},z)$ are convex subsets and
their boundaries are unions of convex subsets in the intersections 
of the inner and outer walls.  
\item By  Proposition \ref{P:structure} (2), if $\pi(\rho)$ is in a chamber then $\rho$ is irreducible.
\end{enumerate}
\end{rmk}


\section{Lagrangian Representations} \label{S:lagrangian}

\subsection{Linear algebra of Lagrangians in $\C^n$}  \label{S:linalg}

 We denote 
by $\Lambda(n)$ the $(n/2)(n+1)$-dimensional manifold of
subspaces of $\C^n$ that are Lagrangian with respect to the standard hermitian structure.
Fixing a preferred Lagrangian
$L_0=\R^n\subset \C^n$, we observe that 
$\Lambda(n)=U(n)/O(n)$, where the orthogonal group $O(n)\subset
U(n)$ is the stabilizer of $L_0$ for the action 
$L_0\mapsto gL_0$.
Define the involution $
\sigma_0(z)\rightarrow \bar{z}
$.
Then to each Lagrangian $L=gL_0=[g]\in \Lambda(n)$ one associates a canonical 
skew-symplectic complex anti-linear involution $\sigma_L:\C^n\rightarrow \C^n$ given by
$
\sigma_L=g\sigma_0 g^{-1}
$, 
whose  set of fixed points is precisely the Lagrangian $L$.
We will set $O_L=$
the stabilizer of $L$, with Lie algebra  ${\mathfrak o}_L$.
Note that $O_L$ is simply the conjugate of $O(n)$ by $g$.
 Let ${\mathfrak u}(n)$ 
denote the Lie algebra of $U(n)$
with the 
$\Ad$-invariant inner product $\langle X,Y\rangle =-\Tr(XY)$.
We have the following useful

\begin{lem}    \label{Ad}
For a Lagrangian $L$:
$
\Ad_{\sigma_L}\bigr|_{{\mathfrak o}_L}={\bf I}
$, and
$
\Ad_{\sigma_L}\bigr|_{{\mathfrak o}_L^\perp}=-{\bf I}
$.
\end{lem}

\Pf
For
 $X\in {\mathfrak u}(n)$, $\Ad_{\sigma_L}(X)$ is  by definition
the derivative at $t=0$ of the curve $\sigma_L e^{tX}\sigma_L\in U(n)$. 
In the case $L=\R^n$, $\sigma_L$ is just complex conjugation, and 
then $\Ad_{\sigma_L}X=\bar X$. 
Using the orthogonal decomposition
$
{\mathfrak u}(n)=i\R^n\oplus {\mathfrak o}(n)\oplus {\mathfrak s}(n)
$,
into diagonal, real orthogonal and symmetric skew-hermitian matrices, 
the result follows immediately.
\EPf

For $g\in U(n)$, let $Z(g)$ denote the centralizer of $g$ with Lie algebra
${\mathfrak z}(g)$.
The relationship between the stabilizers of a pair of Lagrangians
 is given precisely by the following

\begin{prop}\label{decomposition}
Let $L_1$, $L_2$ be two Lagrangian subspaces with stabilizers $O_1$, $O_2$,  and
  let  $g=\sigma_1\sigma_2$ be the composition
of the  corresponding Lagrangian involutions.  Let 
${\mathfrak o}_1$, ${\mathfrak o}_2$ denote the
Lie algebras of $O_1$ and $O_2$. Then
\begin{enumerate}
\item $O_1\cap O_2\subset Z(g)$;
\item There is an orthogonal decomposition
$
{\mathfrak z}(g)=
({\mathfrak o}_1+{\mathfrak o}_2)^\perp \oplus ({\mathfrak o}_1\cap {\mathfrak
o}_2)
$;
\item $2 \dim ({\mathfrak o}_1\cap {\mathfrak o}_2) = \dim{\mathfrak z}(g) -n$.
\end{enumerate}
\end{prop} 
\Pf 
Observe first
 that ${\mathfrak z}(g)=\Ker ({\bf I}-\Ad_g)=
\Ker ({\bf I}-\Ad_{\sigma_1\sigma_2})$.
Using  Lemma \ref{Ad}, we obtain:
$({\mathfrak o}_1+{\mathfrak o}_2)^\perp
 \oplus ({\mathfrak o}_1\cap {\mathfrak
o}_2)\subset{\mathfrak z}(g)
$.
Let $P$ denote the orthogonal projection to ${\mathfrak o}_1\cap
{\mathfrak o}_2$, and let $P_1=(1/2)({\bf I}+\Ad_{\sigma_1})$ and
$P_2=(1/2)({\bf I}+\Ad_{\sigma_2})$ denote the projections to
${\mathfrak o}_1$ and ${\mathfrak o}_2$, respectively.  If $X\in
{\mathfrak z}(g)$, then $\Ad_{\sigma_1}X=\Ad_{\sigma_2}X$, which implies
$P_1X=P_2X$.  Hence, $P\bigr|_{{\mathfrak z}(g)}=
P_1\bigr|_{{\mathfrak z}(g)}=P_2\bigr|_{{\mathfrak z}(g)}$.  In
particular, if $X\in {\mathfrak z}(g)\cap ({\mathfrak o}_1\cap
{\mathfrak o}_2)^\perp$, then $P_1X=P_2X=0$, and $X\in ({\mathfrak
o}_1+ {\mathfrak o}_2)^\perp$.  This proves (2).
Finally, (3) follows from (2).
\EPf

\begin{cor} 
If   $g=\sigma_1\sigma_2$ is regular 
{\rm (}i.e.\ ${\mathfrak z}(g)$ is isomorphic to $i\R^n${\rm )}, then
\begin{enumerate}
\item $O_1\cap O_2= \{ {\bf I}\}$, 
\item 
$O_1\cap Z(g)=O_2\cap Z(g)= \{ {\bf I}\}$.
\end{enumerate}
That is: ${\mathfrak u}(n)=i\R^n\oplus {\mathfrak o}_1\oplus {\mathfrak
o}_2$ (not necessarily orthogonal).
\end{cor}

\begin{dfn} \label{D:tau}
We define three maps:
\begin{align*} 
&\tau_1:\Lambda(n)\longrightarrow U(n): 
L\longmapsto \sigma_L\sigma_0 \ ; \\
&\tau_2: \Lambda^2(n)\longrightarrow U(n): 
(L_1,L_2)\longmapsto\sigma_{L_1}\sigma_{L_2}\ ; \\
&\tau_3 : \Lambda^3(n)\longrightarrow U^2(n) : 
(L_1,L_2,L_3) \longmapsto (\tau_2(L_1,L_2), \tau_2(L_2,L_3))\ . 
\end{align*}
\end{dfn}

\begin{lem}  \label{L:tau}
We have the following:
\begin{enumerate}
\item $ \tau_1([g])= gg^T $; 
\item $\tau_2(L_1,L_2)=
 \tau_1(L_1)\overline{\tau_1(L_2)}$, and $\tau_2(L,L)={\bf I}$; 
\item $\tau_2(L_1,L_3)=\tau_2(L_1,L_2)\tau_2(L_2,L_3)$.
 \end{enumerate}
\end{lem}
We prove some
 elementary facts 
about each of these maps.  Let $S(n)$ denote the space of symmetric $n\times n$ complex 
matrices.

\begin{prop} \label{P:tau1}
The map $\tau_1: \Lambda(n)\rightarrow U(n)$ is an embedding 
with image $U(n)\cap S(n)$.
\end{prop}
\Pf  The fact that the image consists of symmetric matrices is the statement Lemma
\ref{L:tau} (1).
We prove that $\tau_1$ is injective. If $\tau_1([g])=\tau_1([h])$,
then: $gg^T=hh^T$; hence
$
h^{-1}g\in U(n)\cap O(n,\C)
$. 
But $U(n)\cap O(n,\C)=O(n)$, so we conclude that $g\in hO(n)$,
 and $[g]=[h]$.
To prove  $\tau_1$ is an embedding we compute 
its derivative.  
Any variation of  $L$ is determined up to first order
by a variation of the involution $\sigma_L$ of the form
$\sigma_{L(t)}= e^{tX}\sigma_L e^{-tX}$, where $X\in {\mathfrak
u}(n)$.  Then: $\dot\sigma_L=[X,\sigma_L]$, so $\dot\sigma_L\sigma_L\in\im({\bf
I}-\Ad_{\sigma_L})$.  In particular, $\dot\sigma_L\sigma_L=0\iff X\in {\mathfrak o}_L\iff
L(t)\equiv L$.  With this understood, we have
$
\dot\tau_1(L)\tau_1^{-1}(L)=(\dot\sigma_L\sigma_0)(\sigma_0\sigma_L)=\dot\sigma_L\sigma_L
$.
Hence, by the discussion above,  $\tau_1$ is an immersion.
One may show that the image is all of 
$S(n)$ either  by noticing that dimensions agree,  or directly using the
following result, whose proof is straightforward.

\begin{lem} \label{L:sym}
If $g\in U(n)\cap S(n)$ there is $h\in O(n)$ such that $hgh^{-1}$ is diagonal.
\end{lem}
Now take $g$ and $h$ as in the lemma.  Clearly, there exists $k\in U(n)$ such
that $kk^T=hgh^{-1}$.  Then: $\tau_1(hk)=g$.
\EPf

\begin{prop} \label{P:tau2}
 $\tau_2 : \Lambda^2(n)\rightarrow U(n)$ is surjective
and is equivariant with respect to the diagonal action on the domain and 
the conjugation action in the target. 
Over the regular elements of $U(n)$ {\rm (}
i.e.\ those whose eigenvalues have multiplicity one{\rm )}
$\tau_2$ is a fibration  with fiber the torus $T^n$. 
 The general fiber is: 
$\tau^{-1}_2(g)=Z(g)\cap S(n)$, where $Z(g)$ is the centralizer of $g$.
\end{prop}
\Pf
Equivariance is an easy computation.
As a consequence, it suffices to prove the remaining statements for
a diagonal $g\in U(n)$.
For such a $g$ we can solve $g=\tau_2([g_1],[g_2])$, and we may
even assume $g_1$ and $g_2$ are diagonal.  Let $g=h_1h_2$ with 
$h_1=\tau_1([g_1])$ and $h_2=\overline{\tau_1([g_2])}$.
Since $h_2$ is determined by $h_1$ and $\tau_1$ is an embedding,
it suffices to find all possible $h_1$.  Note that since $g$ is diagonal
and $h_1, h_2$ are symmetric, $h_1, h_2\in Z(g)\cap S(n)$.  Conversely,
if $h_1\in  Z(g)\cap S(n)$, then by Proposition \ref{P:tau1},
 $h_1\in \im(\tau_1)$.  Since $h_2=h_1^{-1}g$,
we obtain $h_2^T=g^T({h_1^{-1}})^T=gh_1^{-1}= h_1^{-1}g=h_2$.  We conclude that
$h_2$ is also symmetric, and hence  $h_2\in \im(\tau_1)$.  Thus, $\tau_2^{-1}(g)$
is diffeomorphic to $Z(g)\cap S(n)$.
\EPf 

\noindent
Note that $Z(g)\cap S(n)=S(n_1)\cap U(n_1)\times\dots\times S(n_k)\cap U(n_k)$,
 where $n_i$, for 
$1\leq i\leq k$, are the multiplicities of the eigenvalues of $g$. 
Finally, we determine the image of $\tau_3$.
\begin{dfn} \label{sym}
A pair  $k_1$, $k_2\in U(n)$ is said to be \emph{symmetrizable} if there is 
$g\in U(n)$ such that both  $gk_1g^{-1}$,
$gk_2g^{-1}\in S(n)$.  The set of symmetrizable pairs will be denoted by $\Sym$.
\end{dfn}

\begin{prop} \label{P:tau3}
The image of $\tau_3$ is precisely the set of symmetrizable pairs: $\Sym\subset U^2(n)$.
\end{prop}
\Pf
Clearly if $\tau_3([g_1],[g_2],[g_3])=(h_1,h_2)$, then 
$\tau_3([g_2^{-1}g_1],L_0,[g_2^{-1}g_3])=(g_2^{-1}h_1g_2,g_2^{-1}h_2g_2)$.
But 
$g_2^{-1}h_1g_2=\tau_2([g_2^{-1}g_1], L_0)=\overline{\tau_1([g_2^{-1}g_1])}$
and 
$g_2^{-1}h_2g_2=\tau_2(L_0,[g_2^{-1}g_3])=\overline{\tau_1([g_2^{-1}g_3])}$
which are symmetric.  Therefore $(h_1,h_2)\in \Sym$.
Conversely, suppose  $(h_1,h_2)\in \Sym$, and let $g$ be a matrix 
such that $gh_1g^{-1}, gh_2g^{-1}\in S(n)$.   We can solve
$$
\tau_2([g_1], L_0)={\tau_1([g_1])}=gh_1g^{-1}\ ;\quad
\tau_2(L_0, [g_2])=\overline{\tau_1([g_2])}=gh_2g^{-1}\ .
$$
Then
$\tau_3([g_1],L_0,[g_2])=(gh_1g^{-1},gh_2g^{-1})$.  Since $\tau_3$
is equivariant, acting by
$g^{-1}$ gives the result.  \EPf


\subsection{The space of Lagrangian representations} \label{S:lrep}

We now define the main object of study in this paper.  Fix an integer $\ell\geq 3$. 
 Given the presentation \eqref{gamma}, a representation
$\rho\in\Hom(\Gamma_\ell, U(n))$ is equivalent to a choice of $\ell$ 
matrices whose product is the identity. By Lemma \ref{L:tau} 
(2) and (3),  we therefore have a map
\begin{align}
\tilde\varphi: \Lambda^\ell(n) &\longrightarrow \Hom(\Gamma_\ell, U(n))\ ;
  \label{tildephi} \\
(L_1,\ldots, L_\ell) &\longmapsto (\tau_2(L_1,L_2), \tau_2(L_2, L_3), 
\ldots, \tau_2(L_\ell, L_1))\ .\notag
\end{align}
$U(n)$ acts diagonally on the left of $\Lambda^\ell(n)$, and 
by Proposition \ref{P:tau2}, 
 $\tilde\varphi$ is equivariant
with respect to this action and the left action by conjugation 
of $U(n)$ on $\Hom(\Gamma_\ell, U(n))$.  
Hence, we have an induced map
$$
\varphi: U(n)\backslash\Lambda^\ell(n)\longrightarrow \Rep(\Gamma_\ell, U(n))\ .
$$

Given 
 $\lambda=(L_1, \ldots, L_\ell)
\in \Lambda^\ell(n)$, 
let $Z(\lambda)=O_{L_1}\cap\cdots\cap O_{L_s}\subset U(n)$ denote the stabilizer, and let
${\mathfrak z}(\lambda)$ be its Lie algebra.  Similarly, for 
$\rho\in\Hom(\Gamma_\ell, U(n))$, let $Z(\rho)$ denote its stabilizer with Lie algebra
${\mathfrak z}(\rho)$.  Because of the equivariance of $\tilde\varphi$, $Z(\lambda)\subset
Z(\rho)$, where $\rho=\tilde\varphi(\lambda)$, 
but the two groups are not equal.  For example,
the center $U(1)$ is always
 in $Z(\rho)$ but never in $Z(\lambda)$.  The precise relationship is given by
the following

\begin{lem} \label{inj}
Given $\lambda\in\Lambda^{\ell}(n)$, 
then $\Ker (D\tilde\varphi_\lambda)\subset{\mathfrak u}(n)$,
where ${\mathfrak u}(n)\rightarrow T_\lambda\Lambda^\ell(n)$ via the $U(n)$ action. 
 If $\rho=\tilde\varphi(\lambda)$, then
${\mathfrak z}(\rho)=\Ker (D\tilde\varphi_\lambda)\oplus {\mathfrak z}(\lambda)$.
\end{lem}

\Pf  
Let $\sigma_s=\sigma_{L_s}$, with 
$\sigma_{\ell+1}=\sigma_1$.  Then: $\tilde\varphi(\lambda)=(\gamma_1,\ldots,
\gamma_\ell)$, where $\gamma_s=\sigma_s\sigma_{s+1}$  
(see Definition \ref{D:tau} and \eqref{tildephi}).
Let $\dot\lambda$ be a tangent vector to $\Lambda^\ell(n)$ at $\lambda$.
Expressing the components of the image
$D\tilde\varphi_{\lambda}(\dot\lambda)=(X_1,\ldots, X_s)$ as 
elements of ${\mathfrak u}(n)$, we have: $X_s=\dot\gamma_s\gamma_s^{-1}$.  Hence,
\begin{equation} \label{E:ldef}
X_s=(\dot\sigma_s\sigma_{s+1}+\sigma_s\dot\sigma_{s+1})\sigma_{s+1}\sigma_s
=\dot\sigma_s\sigma_s  +\sigma_s\dot\sigma_{s+1}\sigma_{s+1}\sigma_s\ .
\end{equation}
Since $\sigma_s$ is an involution, we conclude from the equation above that 
$\dot\lambda\in\Ker(D\tilde\varphi_{\lambda})$ if and only if
$\sigma_s\dot\sigma_s=\sigma_{s+1}\dot\sigma_{s+1}$, for all $s=1,\ldots,\ell$.   
As in the proof of Proposition \ref{P:tau1},
 $\sigma_s\dot\sigma_s\in \im({\bf I}-\Ad_{\sigma_s})$. 
If we let $O_s$ denote the stabilizer of the Lagrangian
corresponding to $\sigma_s$, and if ${\mathfrak o}_s$ is the  Lie algebra of $O_s$, then
the kernel of $D\tilde\varphi_\lambda$ is determined by an element in
\begin{align*}
\im({\bf I}-\Ad_{\sigma_1})\cap\cdots\cap \im({\bf I}-\Ad_{\sigma_{\ell}})
&= {\mathfrak o}_1^\perp\cap\cdots\cap {\mathfrak o}_{\ell}^\perp
= ({\mathfrak o}_1+\cdots+ {\mathfrak o}_{\ell})^\perp \\
&= ({\mathfrak o}_1+{\mathfrak o}_2+
{\mathfrak o}_2+{\mathfrak o}_3+\cdots+ {\mathfrak o}_{\ell-1}+
{\mathfrak o}_{\ell})^\perp \\
&= ({\mathfrak o}_1+{\mathfrak o}_2)^\perp\cap\cdots\cap (
{\mathfrak o}_{\ell-1}+{\mathfrak o}_{\ell})^\perp \ .
\end{align*}
By Proposition \ref{decomposition} (2)
 $({\mathfrak o}_s+{\mathfrak o}_{s+1})^\perp\subset {\mathfrak
z}(\gamma_s)$.  
Since
$$
{\mathfrak z}(\rho)= 
{\mathfrak z}(\gamma_1)\cap\cdots\cap {\mathfrak z}(\gamma_{\ell-1})=
({\mathfrak o}_1\cap\cdots\cap {\mathfrak o}_{\ell})\oplus
 ({\mathfrak o}_1+{\mathfrak o}_2)^\perp\cap\cdots\cap (
{\mathfrak o}_{\ell-1}+{\mathfrak o}_{\ell})^\perp \ ,
$$
and ${\mathfrak z}(\lambda)=
{\mathfrak o}_1\cap\cdots\cap {\mathfrak o}_{\ell}$,
the result follows.
\EPf

We take the opportunity to point out a fact about the image of
$D\tilde\varphi_\lambda$.

\begin{lem} \label{L:orth}
Let $(X_1,\ldots, X_\ell)\in \im (D\tilde\varphi_{\lambda})$, with $\lambda$ as above.
Then:  $X_s\in ({\mathfrak o}_s\cap {\mathfrak o}_{s+1})^\perp$ for each
$s=1,\ldots,\ell$.
\end{lem}

\Pf  From Lemma \ref{Ad} and the proof of Lemma \ref{inj}, we have
$$
\dot\sigma_s\sigma_s\in\im({\bf I}-\Ad_{\sigma_s})={\mathfrak o}_s^\perp\ ,\qquad 
\dot\sigma_{s+1}\sigma_{s+1}\in\im({\bf I}-\Ad_{\sigma_{s+1}})=
{\mathfrak o}_{s+1}^\perp
\ .
$$
Now if $Z\in {\mathfrak o}_s\cap {\mathfrak o}_{s+1}$, then by \eqref{E:ldef}
and Lemma \ref{Ad} again,
$$
\langle Z, X_s\rangle=\langle Z, \Ad_{\sigma_s}(\dot\sigma_{s+1}\sigma_{s+1})\rangle
=\langle \Ad_{\sigma_s}Z, \dot\sigma_{s+1}\sigma_{s+1}\rangle
=\langle Z, \dot\sigma_{s+1}\sigma_{s+1}\rangle
=0\ .
$$
\EPf

\begin{dfn} \label{lhom}
A representation $\rho\in  \Hom(\Gamma_\ell, U(n))$ is called a \emph{Lagrangian representation} if it is in the image of $\tilde\varphi$.
We denote the space of \emph{Lagrangian representations} by
$$
\LHom(\Gamma_\ell,U(n))=
\im(\tilde\varphi)
\subset \Hom(\Gamma_\ell, U(n)) \ .
$$
Similarly, 
the image of $\varphi$ is the \emph{moduli space of Lagrangian representations}.
$$
\LRep(\Gamma_\ell,U(n))
=\im(\varphi)
\subset \Rep(\Gamma_\ell, U(n)) \ .
$$
We also set
\begin{align*}
\LHom_{\mathfrak a}(\Gamma_\ell, U(n))&=\LHom(\Gamma_\ell,U(n))\cap  
\Hom_{\mathfrak a}(\Gamma_\ell, U(n))\ ;\\
\LRep_{\mathfrak a}(\Gamma_\ell, U(n))&=\LRep(\Gamma_\ell,U(n))\cap 
 \Rep_{\mathfrak a}(\Gamma_\ell, U(n))\ .
\end{align*}
\end{dfn}

From general considerations of group actions,
 $\Rep^{irr.}(\Gamma_\ell,U(n))$ is a smooth (open)
manifold, since the isotropy $Z(\rho)$ of an irreducible representation $\rho$ is just
 the center of $U(n)$.  Let: $\Lambda^n_{irr.}(n)=\tilde\varphi^{-1}(\Hom^{irr.}(\Gamma_\ell,
U(n))$. Then  
for Lagrangian representations we have the following

\begin{prop} \label{P:embedding}
\begin{enumerate}
\item  For $\lambda\in\Lambda^\ell(n)$ and $\rho=\tilde\varphi(\lambda)$, the fiber
$\tilde\varphi^{-1}(\rho)\simeq Z(\rho)/Z(\lambda)$.  In particular, 
$\LHom^{irr.}(\Gamma_\ell, U(n))$ is an embedded submanifold of dimension
 $$\dim\left(\LHom^{irr.}(\Gamma_\ell,U(n))\right)=
\frac{(\ell-1)}{2}n^2+\frac{\ell}{2} n-1\ ,$$
and: $\tilde\varphi:\Lambda^{\ell}_{irr.}(n)\to
\LHom^{irr.}(\Gamma_\ell, U(n))$ is a circle bundle. 
\item
$U(n)$ acts freely on $\Lambda^n_{irr.}(n)$.
Moreover, 
$$\varphi:U(n)\backslash \Lambda^{\ell}_{irr.}(n)\longrightarrow
\LRep^{irr.}(\Gamma_\ell, U(n))
\subset\Rep^{irr.}(\Gamma_\ell, U(n))
$$
is an embedding with
 $$\dim\left(\LRep^{irr.}(\Gamma_\ell,U(n))\right)=
\frac{(\ell-2)}{2}n^2+\frac{\ell}{2} n\ .$$
\end{enumerate}
\end{prop}

\Pf  We determine the fiber of $\tilde\varphi$.  Suppose
$\rho=\tilde\varphi(\lambda)=\tilde\varphi(\lambda')$, 
where $\lambda=(L_1,\ldots, L_\ell)$ and
 $\lambda'=(L^{\prime}_1,\ldots, L^{\prime}_\ell)$.
 By Propositions
 \ref{P:tau1} and \ref{P:tau2}, $L^{\prime}_1=hL_1$ and $L^{\prime}_2=hL_2$ for
$h\in Z(\rho(\gamma_1))\cap S(n)$. 
 Applying the result to each pair $L_s$, $L_{s+1}$, we see that
in fact
$
h\in Z(\rho(\gamma_1))\cap\cdots\cap Z(\rho(\gamma_{\ell-1}))\cap S(n)
$.
In particular, $h\in Z(\rho)$.  Conversely, by equivariance, 
$Z(\rho)$ acts on the fiber of $\tilde\varphi$ with $Z(\lambda)$.  The remaining statements
follow from Lemma \ref{inj}.
\EPf

We will denote the restriction of 
 the spectral projection  to  the Lagrangian representations
also by 
$\pi: \LHom(\Gamma_\ell, U(n))\to{\mathcal A}^{\Z}_\ell(n)$.
 By analogy with Definition \ref{D:unitary}, we have

\begin{dfn}  \label{D:lagrangian}
Let $\Lellbar^\ast(n)= \pi
(\LHom(\Gamma_\ell, U(n)))\cap \Pellbar^\ast(n)$.
For each collection of multiplicities ${\mathfrak m}=(m^s)$,
 and subsets $z\subset\{1,\ldots, \ell\}$,
we set: $\Lell(n,{\mathfrak m}, z)=\Lellbar^\ast(n)\cap\Pell(n,{\mathfrak m}, z)$.
\end{dfn}

 From the definition we have: $\Lellbar^\ast(n)\subset \Uellbar^\ast(n)$.  
The goal of this paper is to prove that in fact 
$\Lellbar^\ast(n)=\Uellbar^\ast(n)$.  
Assuming Theorem \ref{maintheorem}, however, we may now give the

\medskip
\emph{Proof of Theorem \ref{symmetrizable}.}  
By Theorem \ref{maintheorem}, the conjugacy classes of 
$A_1,\ldots, A_\ell$ may be realized by a Lagrangian representation.
  Hence, we may find $B_i$ as in the statement of 
Theorem  \ref{symmetrizable} such that $B_i=\sigma_{L_i}
\sigma_{L_{i+1}}$ for Lagrangians $L_1,\ldots, L_\ell$,
  where $L_{\ell+1}=L_1$.  In particular, the pair
 $(B_i, B_{i+1})$ is in the image of $\tau_3$ for each $i$.  
The result then follows from Proposition \ref{P:tau3}.
\EPf


\subsection{The symplectic structure} \label{S:symplectic}

The purpose of this section is to show that the tangent space to the Lagrangian
representations for fixed conjugacy classes is isotropic with respect to the
natural symplectic form.  
We begin with a brief review of quasi-Hamiltonian reduction.  For more
details, see \cite{AMM}.
Let $(M,\omega)$ be a manifold equipped with a 2-form $\omega$, $G$ a 
Lie group with Lie algebra
${\mathfrak g}$ and
$
G\times M\rightarrow M
$
a Lie group action preserving $\omega$.  In order to define
a $G$-valued moment map we assume the existence of an Ad-invariant inner product 
$\langle\,  ,\, \rangle$  on ${\mathfrak g}$. 
Let $\theta^R$ and $\theta^L$ be the right and left Maurer-Cartan forms on $G$.
That is, for $V\in T_gG$,  $\theta^L_g(V)=g^{-1}V\in {\mathfrak g}$
 and $\theta^R_g(V)=Vg^{-1}\in {\mathfrak g}$ ($g^{-1}dg$ and $dgg^{-1}$
 in matrix groups).
Let $\chi$ be the bi-invariant closed Cartan 3-form defined by
$$
\chi= \frac{1}{2}\left\langle\theta^L,[\theta^L,\theta^L]\right\rangle=
\frac{1}{2}\left\langle\theta^R,[\theta^R,\theta^R]\right\rangle.
$$

\begin{dfn}
A \emph{quasi-Hamiltonian $G$-space} $(M,G,\omega,\mu)$ is a  manifold equipped with
a 2-form $\omega$ that is invariant under the action of $G$ and an equivariant moment map 
$\mu:M\to G$ satisfying
\begin{enumerate}
\item $d\omega=-\mu^*\chi$
\item $\imath_{\xi^{\#}}\omega =\frac{1}{2}\langle\mu^*(\theta^L+\theta^R),\xi\rangle$
\item $\ker \omega_x=\{\ \xi^{\#}(x)\ |\ \xi\in \ker({\bf I}+\Ad_{\mu(x)})\ \}$.
\end{enumerate}
\end{dfn}
Here, $\xi^{\#}$ denotes the vector field on $M$ induced by
$\xi\in {\mathfrak g}$ and the action of $G$.
The following theorem is proved in \cite{AMM}.

\begin{thm}  \label{T:amm}
Let $(M,G,\omega,\mu)$ be a 
quasi-Hamiltonian space as above.  Let $\imath:\mu^{-1}({\bf I})\rightarrow M$
be the inclusion and $p: \mu^{-1}({\bf I})
\rightarrow M^{red.}= \mu^{-1}({\bf I})/G$  the 
projection on the orbit space.  Then there exists a unique symplectic
form $\omega^{red}$ on the smooth stratum of the reduced space $M^{red}$
such that $p^\ast\omega^{red}= \imath^\ast \omega$ on $\mu^{-1}({\bf I})$.
\end{thm}

This formulation of symplectic reduction is well-adapted to computations on the
representation space of the free group with fixed conjugacy classes.  
Let 
$\Hom_{\mathfrak a}(\Gamma_\ell, U(n))$
and $\Rep_{\mathfrak a}(\Gamma_\ell, U(n))$ be as in Definition \ref{specproj}.
Then $\Hom_{\mathfrak a}(\Gamma_\ell, U(n))$ is 
naturally 
contained in $M_\mathfrak a=C_1 \times \cdots C_\ell$
where $\{C_s\}$ are the conjugacy class
of $U(n)$ prescribed by $\mathfrak a$. Moreover, 
$\Hom_{\mathfrak a}(\Gamma_\ell, U(n))=
\mu^{-1}({\bf I})
$, where $\mu(\gamma_1,\cdots ,\gamma_\ell)=\gamma_1 \gamma_2\cdots
\gamma_\ell
\in U(n)$, and 
$\Rep_{\mathfrak a}(\Gamma_\ell, U(n))=\mu^{-1}({\bf I})/U(n)$.
To describe the form $\omega$, we require

\begin{dfn}
Let $(M_1,\omega_1,\mu_1)$ and $(M_2,\omega_2,\mu_2)$ be two quasi-Hamiltonian
$G$-spaces. Then $M_1\times M_2$ is also a quasi-Hamiltonian $G$-space, called the \emph{fusion product of} $M_1$ and $M_2$.  
The moment map  is given by $\mu_1\mu_2: M_1\times M_2\rightarrow G$, and the  2-form is given by
$
\omega= \omega_1 +\omega_2 +\left\langle\mu^*_1\theta^L\wedge
\mu^*_2\theta^R\right\rangle
$.
\end{dfn}
Explicitly, we have
$$
\left\langle\mu^*_1\theta^L\wedge \mu^*_1\theta^R\right\rangle((v_1,v_2),(w_1,w_2))
= \frac{1}{2}\left( \langle\mu^*_1\theta^L(v_1), \mu^*_2\theta^R(w_2)\rangle
- \langle\mu^*_1\theta^L(w_1), \mu^*_2\theta^R(v_2)\rangle\right ).
$$
To find the expression of the fusion product for a product conjugacy classes,
recall that the fundamental vector field corresponding to $\xi \in {\mathfrak g}$
 at a point $\gamma$ is
$$\xi^{\#}= \xi \gamma -\gamma \xi= ({\bf I}-\Ad_\gamma)\xi \gamma = 
\gamma(\Ad_{{\gamma}^{-1}}-{\bf I})\xi\ .$$
The 2-form on a conjugacy class $C$ is given by
$$
\omega_\gamma(\xi^{\#},\eta^{\#})=\frac{1}{2}\left
 ( \langle \Ad_\gamma \xi, \eta\rangle-\langle \Ad_\gamma \eta, \xi\rangle \right ).
$$
For the product 
of two conjugacy classes $C_1$ and $C_2$, let
$\mu_i: C_i\rightarrow G$ be the tautological embeddings.  Then
\begin{align*}
\mu_1^*\theta^L(\xi_1^{\#})&=
\theta^L({\mu_1}_*\xi_1^{\#})=\theta^L(\xi_1^{\#})=
\theta^L(\gamma_1(\Ad_{\gamma_1^{-1}}-{\bf I})\xi_1)\\
&=
\gamma_1^{-1}\gamma_1(\Ad_{\gamma_1^{-1}}-{\bf I})
\xi_1=(\Ad_{\gamma_1^{-1}}-{\bf I})\xi_1\ .
\end{align*}
Similarly, $\mu_2^*\theta^R(\eta_2^{\#})=
({\bf I}-\Ad_{\gamma_2})\eta_2$.
Using these formulas, the 2-form 
 on the product $C_1\times C_2$ of two conjugacy classes is
\begin{align*}
\omega_{(\gamma_1,\gamma_2)}&\left((\xi_1^{\#},\xi_2^{\#}),(\eta_1^{\#},
\eta_2^{\#}) \right)=
\frac{1}{2}\left
 ( \langle \Ad_{\gamma_1} \xi_1, \eta_1\rangle -\langle \Ad_{\gamma_1} \eta_1,
 \xi_1\rangle \right ) \\
& + \frac{1}{2}\left
 ( \langle \Ad_{\gamma_2} \xi_2, \eta_2\rangle-\langle \Ad_{\gamma_2} \eta_2,
 \xi_2\rangle \right ) +\frac{1}{2}\langle 
({\bf I}-\Ad_{\gamma_1})\xi_1,\Ad_{\gamma_1}({\bf I}-\Ad_{\gamma_2})\eta_2\rangle
 - \{\xi\leftrightarrow\eta\}
\end{align*}
where $\xi\leftrightarrow\eta$ means that the previous terms are repeated with
$\xi$ and $\eta$ interchanged,
 keeping the indices unchanged.
In general, for the product $C_1\times\cdots\times C_\ell$ we obtain
\begin{align*}
&\omega_{(\gamma_1,\cdots ,\gamma_\ell)}\bigl((\xi_1^{\#},\cdots, 
\xi_\ell^{\#}),(\eta_1^{\#},\cdots ,\eta_\ell^{\#})\bigr)=\\
&\qquad=\frac{1}{2}
 \biggl\{ \sum _{s=0}^{\ell} \langle \Ad_{\gamma_s} \xi_s, \eta_s\rangle 
+  \sum_{t=1}^{\ell-1}
\bigl\langle ({\bf I}-\Ad_{\gamma_1})\xi_1+ \Ad_{\gamma_1}({\bf I}-\Ad_{\gamma_2})
\xi_2 +\cdots \\
&\qquad\qquad\cdots 
+\Ad_{\gamma_1\cdots \gamma_{t-1}}({\bf I}-\Ad_{\gamma_{t}})\xi_{t}, 
\Ad_{\gamma_1\cdots \gamma_{t}}({\bf I}-\Ad_{\gamma_{t+1}})\eta_{t+1}
\bigr\rangle \biggr\}-\{\xi\leftrightarrow\eta\} \\ 
&\qquad= \frac{1}{2}
 \biggl\{ \sum _{s=0}^{\ell} \langle \Ad_{\gamma_s} \xi_s, \eta_s\rangle+
+\sum_{0\leq s<t\leq \ell-1}\left \langle 
 \Ad_{\gamma_1\cdots \gamma_{s}}({\bf I}-\Ad_{\gamma_{s+1}})\xi_{s+1}, 
\Ad_{\gamma_1\cdots \gamma_{t}}({\bf I}-\Ad_{\gamma_{t+1}})\eta_{t+1}
\right\rangle\biggr\}\\
&\qquad\qquad -\{\xi\leftrightarrow\eta\}\ .
\end{align*}
\begin{prop} \label{P:amm}
The product of conjugacy classes of a compact Lie group $G$,
 $C_1\times\cdots\times C_\ell$ is a quasi-Hamiltonian
space equipped with the moment map which is the product of the embeddings in $G$ and
the following 2-form: 
\begin{align*}
\omega_{(\gamma_1,\cdots ,\gamma_\ell)}&\bigl( (\xi_1^{\#},\cdots,
\xi_\ell^{\#}),
(\eta_1^{\#},\cdots ,\eta_\ell^{\#})\bigr )=
 \frac{1}{2}
 \biggl \{ \sum _{s=0}^{\ell} ( \Ad_{\gamma_s} \xi_s, \eta_s) +\\
& +\sum_{0\leq s<t\leq\ell-1}\left ( \Ad_{\gamma_1\cdots \gamma_{s}}
 ({\bf I}-\Ad_{\gamma_{s+1}})\xi_{s+1}, 
\Ad_{\gamma_1\cdots \gamma_{t}}({\bf I}-\Ad_{\gamma_{t+1}})\eta_{t+1}
\right )\biggr \}
-\{\xi\leftrightarrow\eta\}\ .
\end{align*}
\end{prop}

\begin{prop} \label{P:isotropic}
The moduli space of moduli of  Lagrangian representations 
$$\LRep_{\mathfrak a}(\Gamma_\ell, U(n))\subset
 \Rep_{\mathfrak a}(\Gamma_\ell, U(n))\ $$ is  isotropic with respect to the
 symplectic structure defined by Proposition \ref{P:amm} and Theorem
 \ref{T:amm}.
\end{prop}
\Pf Let
$$
 X_s=\dot\sigma_s\sigma_s+\Ad_{\sigma_s}(\dot\sigma_{s+1}\sigma_{s+1})\
 ,\qquad
 Y_s=\dot\rho_s\rho_s+\Ad_{\sigma_s}(\dot\rho_{s+1}\rho_{s+1})\ .
$$
where $\rho_s=\sigma_s$ (see \eqref{E:ldef}). 
By the assumption of fixed conjugacy classes, we have
\begin{alignat*}{2} 
\dot\sigma_s\sigma_s &= \xi_s-\Ad_{\sigma_s}\xi_s\ , &\qquad
\dot\rho_s\rho_s &= \eta_s-\Ad_{\sigma_s}\eta_s\\
\dot\sigma_{s+1}\sigma_{s+1} &= \xi_s-\Ad_{\sigma_{s+1}}\xi_s\ , &\qquad 
\dot\rho_{s+1}\rho_{s+1} &= \eta_s-\Ad_{\sigma_{s+1}}\eta_s\ .
\end{alignat*}
In particular,
\begin{equation} \label{E:adx}
\Ad_{\sigma_s} X_s=\dot\sigma_{s+1}\sigma_{s+1}-\dot\sigma_{s}\sigma_{s}
\ ,\quad \Ad_{\sigma_s} Y_s=\dot\rho_{s+1}\rho_{s+1}-\dot\rho_{s}\rho_{s}\ .
\end{equation}
It follows that
\begin{align*}
\langle \Ad_{\gamma_s}\xi_s,\eta_s\rangle &= \langle
\Ad_{\sigma_{s+1}}\xi_s, \Ad_{\sigma_s}
\eta_s\rangle
=\langle \xi_s-\dot\sigma_{s+1}\sigma_{s+1}, \eta_s-\dot\rho_s\rho_s\rangle\\
&=\langle\xi_s,\eta_s\rangle+\langle \dot\sigma_{s+1}\sigma_{s+1},\dot\rho_s\rho_s\rangle
-\langle \xi_s,\dot\rho_s\rho_s\rangle -\langle \eta_s,\dot\sigma_{s+1}\sigma_{s+1}\rangle\ .
\end{align*}
Notice that since $\dot\rho_s\rho_s$ is in the ($-1$)-eigenspace of $\Ad_{\sigma_s}$,
$$
2\langle \xi_s, \dot\rho_s\rho_s\rangle = \langle \xi_s-\Ad_{\sigma_s}\xi_s, \dot\rho_s\rho_s\rangle=\langle \dot\sigma_s\sigma_s, \dot\rho_s\rho_s\rangle\ .
$$
Similarly, $ 2\langle \eta_s,\dot\sigma_{s+1}\sigma_{s+1}\rangle= \langle \dot\rho_{s+1}\rho_{s+1},\dot\sigma_{s+1}\sigma_{s+1}\rangle$.  Because of the symmetry upon interchanging $\sigma$ and $\rho$, these terms cancel, and  we are left with
\begin{equation} \label{E:term1}
\sum_{s=1}^\ell \langle \Ad_{\gamma_s}\xi_s,\eta_s\rangle-\langle \Ad_{\gamma_s}\eta_s,\xi_s\rangle
=
\sum_{s=1}^\ell \langle \dot\sigma_{s+1}\sigma_{s+1},\dot\rho_s\rho_s\rangle-\langle \dot\sigma_s\sigma_s,\dot\rho_{s+1}\rho_{s+1}\rangle\ .
\end{equation}
For the second term, notice that for a Lagrangian representation
$\gamma_1\cdots\gamma_s=\sigma_1\sigma_{s+1}$.  Hence,
\begin{align*}
\sum_{0\leq s<t\leq \ell-1}
\langle \Ad_{\gamma_1\cdots\gamma_s}X_{s+1},\Ad_{\gamma_1\cdots\gamma_{t}}Y_{t+1}\rangle
&=\sum_{0\leq s<t\leq \ell-1}
\langle \Ad_{\sigma_{s+1}}X_{s+1},\Ad_{\sigma_{t+1}}Y_{t+1}\rangle\\
&=\sum_{1\leq s<t\leq \ell}
\langle \Ad_{\sigma_{s}}X_{s},\Ad_{\sigma_{t}}Y_{t}\rangle\ .
\end{align*}
Using \eqref{E:adx} (and recalling the convention that $\rho_{\ell+1}=\rho_1$) we have
\begin{align*}
\sum_{0\leq s<t\leq \ell-1}\langle \Ad_{\gamma_1\cdots\gamma_s}&X_{s+1},
\Ad_{\gamma_1\cdots\gamma_{t}}Y_{t+1}\rangle =
\sum_{1\leq s<t\leq \ell}\langle\dot\sigma_{s+1}\sigma_{s+1}-\dot\sigma_s\sigma_s,
\dot\rho_{t+1}\rho_{t+1}-\dot\rho_t\rho_t\rangle\\
&=\sum_{1\leq s\leq \ell-1}\langle\dot\sigma_{s+1}\sigma_{s+1}-\dot\sigma_s\sigma_s,
\dot\rho_{1}\rho_{1}-\dot\rho_{s+1}\rho_{s+1}\rangle\\
&=\sum_{1\leq s\leq \ell-1}  \langle\dot\sigma_s\sigma_s,
\dot\rho_{s+1}\rho_{s+1}\rangle-\langle\dot\sigma_{s+1}\sigma_{s+1},
\dot\rho_{s+1}\rho_{s+1}\rangle+\langle\dot\sigma_{s+1}\sigma_{s+1}-\dot\sigma_s\sigma_s,\dot\rho_1\rho_1\rangle\\
&=\langle\dot\sigma_{\ell}\sigma_{\ell}-\dot\sigma_1\sigma_1,\dot\rho_1\rho_1\rangle+\sum_{1\leq s\leq \ell-1}  \langle\dot\sigma_s\sigma_s,
\dot\rho_{s+1}\rho_{s+1}\rangle-\langle\dot\sigma_{s+1}\sigma_{s+1},
\dot\rho_{s+1}\rho_{s+1}\rangle\\
&=\sum_{1\leq s\leq \ell}     \langle\dot\sigma_s\sigma_s,
\dot\rho_{s+1}\rho_{s+1}\rangle-\langle\dot\sigma_s\sigma_s,
\dot\rho_{s}\rho_{s}\rangle
\end{align*}
Hence,
\begin{align*}
\sum_{0\leq s<t\leq\ell-1}\langle
\Ad_{\gamma_1\cdots\gamma_s}X_{s+1},\Ad_{\gamma_1\cdots\gamma_{t}}Y_{t+1}\rangle
&-\langle \Ad_{\gamma_1\cdots\gamma_s}Y_{s+1},\Ad_{\gamma_1\cdots\gamma_{t}}X_{t+1}\rangle\\
&=\sum_{s=1}^\ell\langle\dot\sigma_s\sigma_s,
\dot\rho_{s+1}\rho_{s+1}\rangle-\langle\dot\sigma_{s+1}\sigma_{s+1},\dot\rho_s\rho_s
\rangle\ .
\end{align*}
The proposition now  follows by comparing this with \eqref{E:term1}.
\EPf


\subsection{The Maslov index} \label{maslov}

In this section, we briefly digress to explain the relationship between the quantity
$I(\rho)$, which we have called the index of a representation, and the usual Maslov index
of a triple of Lagrangians, in the case $\rho$ is a Lagrangian
representation.
The diagonal action of the symplectic group acting on triple 
of Lagrangian subspaces $(L_1,L_2,L_3)$  in $\C^n$
 has a finite number of orbits. 
 To classify the orbits, one introduces
the notion of an \emph{inertia index} (or \emph{Maslov index}) of 
a Lagrangian triple (cf.\ \cite[p.\ 486]{KS}).

\begin{dfn} The \emph{inertia index} $\tau(\lambda)$ of a triple 
${\lambda}=(L_1, L_2, L_3)$ of 
Lagrangian subspaces of $\C^n$ is the signature of the quadratic
form $q$ defined on the $3n$ (real)
 dimensional vector space $L_1 \oplus L_2 \oplus
L_3$ by: $q(x_1,x_2,x_3)=\omega(x_1,x_2)+\omega(x_2,x_3)+\omega(x_3,x_1)$, 
where $\omega$ is the standard symplectic form on $\C^n$.
\end{dfn}

\noindent In order to state the
 symplectic classification of triples of Lagrangians,
 we need the following data.
For ${\mathfrak d}=(n_0,n_{12},n_{23},n_{31},\tau)\in \N^4\times \Z$, let
$ C_{\mathfrak d} $ denote the set of all   
${\lambda}=(L_1,L_2,L_3)$ satisfying
$\tau({\lambda})=\tau$, 
$\dim(L_1 \cap L_2 \cap L_3) = n_0$, and
$\dim(L_j\cap L_k)
=n_{jk} $.  For the following result, see \cite[p.\
493]{KS}.
\begin{prop}
\label{sympclass}
$C_{\mathfrak d}$ is non-empty if and only if ${\mathfrak d}=
(n_0,n_{12},n_{23},n_{31},\tau)$ satisfies the
conditions
\begin{enumerate}
\item $0\leq n_0\leq n_{12},n_{23},n_{31}\leq n $.
\item $n_{12}+n_{23}+n_{31}\leq n+2n_0$.
\item $|\tau|\leq n+2n_0-(n_{12}+n_{23}+n_{31})$.
\item $\tau\equiv n-(n_{12}+n_{23}+n_{31}) \mod {2 \Z}$.
\end{enumerate}
If $\lambda$ and $\lambda'$ are two triples of Lagrangian subspaces of $\C^n$, there
exists a symplectic map $\psi\in Sp(\C^n)$ such that $\psi(L_1)=L'_1$,
$\psi(L_2)=L'_2$ and $\psi(L_3)=L'_3$, if and only if $n_0=n'_0$, $n_{12}=n'_{12}$,
$n_{23}=n'_{23}$, $n_{31}=n'_{31}$ and $\tau=\tau'$.
\end{prop}

Using this classification one may show

\begin{prop}[{\cite[Theorem 4.4]{FMS}}]    \label{indcomp}
Let $\lambda=(L_1,L_2,L_3)$, 
$\rho=\tilde\varphi(\lambda)$, and 
$n_{jk}=\dim(L_j\cap L_k)$.
Then
 $$\tau(\lambda)=3n-2 I(\rho)-(n_{12}+n_{23}+n_{31})\ .$$
\end{prop}

This relationship between $\tau$ and $I$ 
gives an alternative proof of Theorem \ref{indexbounds1} for the case $\ell=3$ (and
assuming Theorem \ref{maintheorem}).

\begin{cor} \label{indexbounds2} 
Let $\lambda$ be a triple of Lagrangian subspaces of $\C^n$, 
$\rho=\tilde\varphi(\lambda)$. Then
$$
n-N_0(\rho)\leq I(\rho)\leq 2n + N_0(\rho) - N_1(\rho)\ .
$$
\end{cor}
\Pf This follows from 
Propositions \ref{sympclass} and \ref{indcomp}, 
and the fact that $N_0(\rho)=n_0$, and $N_1(\rho)= 
n_{12}+n_{23}+n_{31}$.
\EPf

The Maslov index generalizes to multiple Lagrangians as follows.
Let $L_1,\ldots , L_\ell$, $\ell\geq 3$, 
be a collection of Lagrangian subspaces in $\C^n$.  We define
$$
\tau(L_1,\ldots , L_\ell)=\tau(L_1,L_2,L_3)+\tau(L_1,L_3,L_4)+
\cdots +\tau(L_1,L_{\ell-1},L_\ell)
\ .
$$
For the next result, set
$
I(L_1,\ldots , L_\ell)= I(\tilde\varphi(L_1,\ldots , L_\ell))
$.

\begin{prop}
Let $L_1,\ldots , L_\ell$, $\ell\geq 4$, be a collection of Lagrangian subspaces
in $\C^n$.  Write $n_{1i}=\dim(L_1\cap L_i)$,  then
$$
I(L_1,\ldots , L_\ell)=I(L_1,L_2,L_3)+I(L_1,L_3,L_4)+\cdots +
I(L_1,L_{\ell-1},L_\ell)-\sum_{i=3}^{\ell-1} (n-n_{1i})\ .
$$
\end{prop}
\Pf 
Observe that 
 if $\spec(\sigma_{L_1}\sigma_{L_3})=
(0,\ldots,0,\alpha_{n_{13}+1},\ldots,\alpha_n)$ then
$$ \spec(\sigma_{L_3}\sigma_{L_1})=
 (0,\ldots,0,
1-\alpha_n, \ldots,
1-\alpha_{n_{13}+1})
\ .
$$
 Summing all the angles in both spectra 
gives us: $n-n_{13}$.  This implies that
$$
I(L_1,L_2,L_3,L_4)=I(L_1,L_2,L_3)+I(L_1,L_3,L_4)-(n-n_{13})\ .
$$
The general case follows by induction.
\EPf

A relationship between $\tau$ and $I$ still exists.  Indeed, this
follows directly from the previous result and  Proposition \ref{indcomp}.

\begin{prop}  \label{iandtau}
For $\ell\geq 3$,
$
\tau(L_1,\ldots , L_\ell)=n\ell-2I(L_1,\ldots , L_\ell)-
(n_{12}+n_{23}+\cdots +n_{\ell 1})
$.
\end{prop}

It is not immediately clear how to prove the  analogue of
Proposition \ref{sympclass} for $\ell\geq 4$, 
since the invariants no longer necessarily classify 
$\ell$-tuples of Lagrangians.
On the other hand, we can \emph{use}  Theorem \ref{indexbounds1}, along with
Proposition \ref{iandtau}, to
prove bounds on the generalized Maslov index.

\begin{thm} \label{T:maslov}
For any $\ell$-tuple of Lagrangians,
$$
|\tau(L_1,\ldots, L_\ell)|\leq  n(\ell-2)+2n_0-(n_{12}+n_{23}+\cdots+n_{\ell
1})\ .
$$
\end{thm}


\section{Deformations of Unitary and Lagrangian Representations} \label{S:deformation}


\subsection{The deformation space}

For an algebraic group $G$ and a finitely presented group $\Gamma$, 
let $\Hom(\Gamma,G)$ be the space homomorphisms
 of $\Gamma$ into $G$. 
If $\Gamma$ has generators $\{\gamma_1,\ldots, \gamma_\ell\}$,
then
$\Hom(\Gamma,G)$ is given the structure of an algebraic variety
as the common locus of 
 inverse images of the identity in $G^\ell$
for a finite number of functions $r_i:G^\ell\rightarrow G$.
  The tangent space to $G^\ell$ 
is identified 
with ${\mathfrak g}^\ell$, where $\mathfrak g$ is the Lie algebra of $G$, 
 by right invariant vector fields.
 If $\rho_t$ is a path of representations,
$\rho_0=\rho$, then
differentiating $\rho_t$ on a word $\gamma_{i_1}\cdots\gamma_{i_m}$, and using
$ X_k= \dot\rho_0(\gamma_{i_k}) \rho_0^{-1}(\gamma_{i_k})$,  we obtain the cocycle
relation
$$
X_1+\Ad_{\rho(\gamma_{i_1})}X_2+\cdots +\Ad_{\rho(\gamma_{i_1}\cdots
\gamma_{i_{m-1}})}X_m=0.
$$
This formula implies the following observation of Weil \cite{W}.
\begin{prop}
The Zariski tangent space $T_\rho \Hom(\Gamma,G)$ is isomorphic to
$Z^1(\Gamma,\Lg)$.
\end{prop}
In order to analyze deformations  fixing conjugacy classes 
we compute the derivative of the curve
 $t\rightarrow \Ad_{e^{tX}}\rho(\gamma)$
 to obtain $\{X-\Ad_{\rho(\gamma)}X\}\rho(\gamma)$.  Identifying 
this with   $X-\Ad_{\rho(\gamma)}X\in \Lg$, we obtain a boundary in the
group cohomology.

We apply these general considerations to the case of $U(n)$ representations of 
 the  free groups $\Gamma=\Gamma_{\ell}$ with presentation as in \eqref{gamma}.
 $\Hom(\Gamma_\ell, U(n))$ is a smooth manifold 
of dimension $(\ell-1)n^2$, 
with tangent space
  at a representation $\rho$
given by
\begin{equation} \label{tangentspace1}
X_1+\Ad_{\rho(\gamma_1)}X_2+\cdots
+\Ad_{\rho(\gamma_1\cdots\gamma_{\ell-1})}X_\ell=0.
\end{equation}



We will be concerned with $g\in U(n)$  with 
fixed multiplicity for the
eigenvalues.  
As in Section \ref{S:lagrangian},
 let ${\mathfrak z}(g)$ be the Lie algebra of the centralizer of $g$.
Then ${\mathfrak z}(g)={\mathfrak u}(\mu_1)\times\cdots\times{\mathfrak u}(\mu_l)$,
where: $\mu_1,\ldots, \mu_l$ are the multiplicities of the eigenvalues of $g$.
We define ${\mathfrak z}^{ab.}(g)\subset{\mathfrak z}(g)
={\mathfrak u}(1)\times\cdots\times{\mathfrak u}(1)$,
 to be the subalgebra
consisting of elements that are block diagonal with respect to this
decomposition.  Alternatively, it is the maximal abelian ideal of ${\mathfrak
z}(g)$.
We have the following

\begin{lem} \label{L:multiplicity}
Let $m$ be a multiplicity structure as in Section
\ref{S:conjugacy}.  Let $U(n,m)$ denote the set of all $g\in U(n)$ with
multiplicity structure $m$.  Then $U(n,m)$ is a smooth submanifold with
tangent bundle {\rm (}identified with a subspace of ${\mathfrak u}(n)${\rm )}
 given by: 
 ${\mathfrak u}(n,m)={\mathfrak z}^{ab.}(g)\oplus {\mathfrak z}(g)^\perp$.
Similarly, if $U(n,m,0)$ is the set of all $g\in U(n,m)$ with
$0\in \spec(g)$, then $U(n,m,0)$ is a smooth submanifold with
tangent bundle  given by
 ${\mathfrak u}(n,m,0)={\mathfrak z}^{ab.,0}(g)\oplus {\mathfrak z}(g)^\perp$,
where the superscript indicates that the first ${\mathfrak u}(1)$ factor is
zero.
\end{lem}

\Pf  It suffices the prove the statement concerning the tangent space.
But small deformations of the eigenvalues are obtained by $g(t)=e^{tX}g$ for
$X\in {\mathfrak z}^{ab.}(g)$.  Conjugating by an arbitary unitary matrix, we find
$$
{\mathfrak u}(n,m)=
\left\{X+({\bf I}-\Ad_g)Y :
X\in  {\mathfrak z}^{ab.}(g)\ ,\ Y\in {\mathfrak u}(n)\right\}\ .
$$
Since ${\mathfrak z}(g)^\perp=\im({\bf I}-\Ad_g)$, the result follows.
The reasoning for $U(n,m,0)$ is similar.
\EPf

Now we prove

\begin{prop}[{cf.\ \cite[Section 5]{MS}}]  \label{maximalrank1} 
Let $\rho:\Gamma_\ell\to U(n)$ be irreducible with
$\pi(\rho)={\mathfrak a}\in{\mathcal U}_{I,\ell}(n,{\mathfrak m}, z)$.
Then near $\rho$, 
$\Rep^{irr.}_{\mathfrak a}(\Gamma_\ell, U(n))$
is a smooth manifold of dimension
$$
\dim\left(\Rep^{irr.}_{\mathfrak a}(\Gamma_\ell, U(n))\right)=
(\ell-2)n^2+2-\sum_{s=1}^\ell\sum_{j=1}^{l_s} (\mu^s_j)^2
$$
 Here, $\mu^s_j$ denotes the
 multiplicity $m^s_j-m^s_{j-1}$ of the $j$-th distinct eigenvalue of $\rho(\gamma_s)$,
$j=1,\ldots, l_s$,  (see Section \ref{S:conjugacy}).
Moreover, the spectral projection
$$
\pi: \Rep^{irr.}(\Gamma_\ell, U(n))\cap \pi^{-1}\left( 
{\mathcal U}_{I,\ell}(n,{\mathfrak m}, z)\right)
\longrightarrow
{\mathcal P}_{I,\ell}(n,{\mathfrak m}, z)\ ,
$$
is locally surjective and is a fibration near $\rho$.
\end{prop}

\Pf
We fix the conjugacy classes of $\rho(\gamma_s)$ for $s\geq 2$ and
determine the variation in $\rho(\gamma_1)$. The space
$\Hom_{{\mathfrak a}^{(1)}}(\Gamma_\ell,U(n))$ of
representations $\rho'$,  $ \rho'(\gamma_s)\simeq
g_s$ for $s\geq 2$, is clearly a manifold of dimension
\begin{equation}  \label{homdim}
\dim\left(\Hom_{{\mathfrak a}^{(1)}}(\Gamma_\ell,U(n))\right)
=\sum_{s=2}^\ell \dim(U(n)/Z(g_s))=
(\ell-1)n^2- 
\sum_{j=2}^{l_s} (\mu^s_j)^2\ ,
\end{equation}
where we have used that: $ \dim Z(g_s)=
\sum_{j=1}^{l_s}(\mu^s_j)^2$.
 We compute the derivative of the map
$$
\pi^{(1)}:\Hom_{{\mathfrak a}^{(1)}}(\Gamma_\ell,U(n))\longrightarrow U(n):
\rho\mapsto\rho(\gamma_1)\ .
$$
Note that $\pi^{(1)}$ takes values in a single fiber of the determinant map.
By \eqref{tangentspace1},  the tangent space to 
$\Hom_{{\mathfrak a}^{(1)}}(\Gamma_\ell,U(n))$ at
 $\rho$  is given by $(X_1,\cdots, X_\ell)\in \mathfrak{g}^\ell$ satisfying
the conditions
$
X_s\in \im({\bf I}-\Ad_{\rho(\gamma_s)})
$, for $s\geq 2$,  and
$$
 X_1\in V^{(1)}= \Ad_{\rho(\gamma_1)}\im({\bf I}-\Ad_{\rho(\gamma_2)})+\cdots
+ \Ad_{\rho(\gamma_1\gamma_2\cdots\gamma_{\ell-1})}\im({\bf
I}-\Ad_{\rho(\gamma_\ell)}) \ .
$$
We claim that $V^{(1)}= z(\rho)^{\perp}$.
Indeed,
\begin{align*}
(V^{(1)})^\perp &=\left\{ \Ad_{\rho(\gamma_1)}\im({\bf
I}-\Ad_{\rho(\gamma_2)})+\cdots +
\Ad_{\rho(\gamma_2\cdots\gamma_{\ell-1})}\im({\bf
I}-\Ad_{\rho(\gamma_\ell)})\right\}^\perp\\
&= \Ad_{\rho(\gamma_1)}\bigl\{(\im({\bf I}-\Ad_{\rho(\gamma_2)}))^\perp\cap
\Ad_{\rho(\gamma_2)}(\im({\bf I}-\Ad_{\rho(\gamma_3)}))^\perp\cdots \\
&\hskip2in\cdots\cap
 \Ad_{\rho(\gamma_2\cdots\gamma_{\ell-1})}(\im({\bf
I}-\Ad_{\rho(\gamma_\ell)}))^\perp\bigr\}\ . 
\end{align*}
Now $\im({\bf I}-\Ad_{\rho(\gamma_2)})^\perp=\mathfrak{z}(\gamma_2)$, and
therefore
$$
\bigl(\im({\bf I}-\Ad_{\rho(\gamma_2)})\bigr)^\perp\cap
\Ad_{\rho(\gamma_2)}\bigl(\im({\bf
I}-\Ad_{\rho(\gamma_3)})\bigr)^\perp
=\mathfrak{z}(\gamma_2)\cap \mathfrak{z}(\gamma_3)\ .
$$
Continuing in this way, we find
$
(V^{(1)})^\perp = \mathfrak{z}(\gamma_2)\cap\cdots\cap
\mathfrak{z}(\gamma_\ell)=\mathfrak{z}(\rho)
$.
We have shown that
$\im(D\pi^{(1)}_\rho)=\mathfrak{z}(\rho)^\perp$.
Hence, if the representation is irreducible (i.e.\ if ${\mathfrak
z}(\rho)=
{\mathfrak z}(U(n))\simeq i\R$), then  by transversality we conclude that
$\Hom_{\mathfrak a}(\Gamma_\ell,U(n))$ is a manifold at an irreducible.
Transversality  applied to 
the product over $s$ of $U(n,m^s)$ (or $U(n,m^s,0)$ if $s\in z$)
also gives 
the statement about local surjectivity and the fibration
structure over the multiplicity space (see Lemma \ref{L:multiplicity}).
For the dimension, we observe that
\begin{itemize}
\item
$ \dim{\mathfrak z}(\rho(\gamma_1))=
\sum_{j=1}^{l_1}(\mu^1_j)^2$; 
\item since $Z(\rho)=Z(U(n))=U(1)$ by irreducibility,
$n^2-1$ is the dimension of $U(n)$-orbit through $\rho$.
\end{itemize}
The dimension of $\pi^{-1}({\mathfrak a})$ is computed by 
subtracting these from \eqref{homdim}.
Since this dimension depends only on the multiplicity structure, it is
constant over the fixed multiplicity space; hence, the
smoothness.  This completes the proof.
\EPf

\begin{rmk}  \label{ms}
The surjectivity in Proposition
 \ref{maximalrank1} also follows from the Mehta-Seshadri Theorem \cite{MS}
which describes irreducible representations with fixed conjugacy classes
in terms of stable parabolic vector bundles.
In the next section we will see that a similar result holds even if we restrict
$\pi$ to the Lagrangian representations,
where we apparently have no such holomorphic description.
\end{rmk}


\subsection{Twisting and bending deformations of Lagrangian representations}

We approach the deformation theory of Lagrangian representations
 by introducing two special families:
twist deformations and real bendings.  Twist deformations are rather simple and apply
equally well to unitary representations, while the bending deformations are particular to
Lagrangian representations.

\begin{dfn} \label{D:twist}
Let $\lambda=(L_1,\ldots, L_\ell)\in\Lambda^\ell(n)$, and $\rho=\tilde\varphi(\lambda)$. 
A \emph{twist deformation} of the Lagrangian representation $\rho$ is a Lagrangian
representation of the 
form: $\rho_\tau=\tilde\varphi(\lambda_\tau)$, where
$\lambda_\tau=(\tau_1 L_1,\ldots, \tau_\ell L_\ell)$ for  some
$\tau=(\tau_1,\ldots, \tau_\ell)\in U^{\ell}(1)$.
\end{dfn}

\begin{rmk} \label{R:det}
Since $\tilde\varphi$ always has the center of $U(n)$ as a fiber, the twist deformations
naturally depend on $\ell-1$ parameters in $U(1)$.
\end{rmk}

The following result is a calculation using the method in the proof of  Lemma \ref{inj}.

\begin{lem} \label{L:twist}
Let ${\mathcal T}_\rho\subset T_\rho\LHom(\Gamma_\ell, U(n))$ denote the subspace tangent
to the twist deformations of $\rho$.
Then
$
{\mathcal T}_\rho=\left[ {\mathfrak u}(1)\times\cdots\times {\mathfrak u}(1)\right]_0
$,
where ${\mathfrak u}(1)$ is the Lie algebra of the center $U(1)\subset U(n)$, and the
subscript $0$ indicates that the sum of the entries vanishes.
\end{lem}

\begin{dfn} \label{D:bending}
Let $\lambda$, $\rho$ be as in Definition \ref{D:twist}.  A \emph{real bending}
 of the Lagrangian representation $\rho$ is a Lagrangian representation of the 
form $\rho_b=\tilde\varphi(\lambda_b)$, where
$$\lambda_b=(L_1,\ldots,L_s, bL_{s+1}, \ldots, bL_{s+r}, L_{s+r+1},\ldots,   L_\ell)$$
 for  some $s, r=1,\ldots, \ell$, and 
$b\in O_{L_s}$ (as usual, we reduce  mod $\ell$ any index greater than $\ell$). 
 Given $s,r$ we shall
say the bending is about $L_s$ and has length $r$.
\end{dfn}

 The twist deformations considered above are special cases
of the action of the group $\Hom(\Gamma, Z(G))$ on 
$\Hom(\Gamma,G)$, and they were
considered in \cite{LM}.
Bending deformations  are inspired by  generalizations of Fenchel-Nielsen twists
  defined by Thurston (see
 \cite{G2} and \cite{JM}).  An important difference is that  the  bending deformations defined in these references \emph{fix} the conjugacy classes of $\rho(\gamma_s)$, whereas those  in Definition \ref{D:bending} change
certain conjugacy classes in a controlled way.

Indeed, from the definition we see 
that a real bending of length $r$ about $L_s$ has the form
\begin{equation} \label{E:bending}
\rho_b(\gamma_{s'})=\begin{cases}
\rho(\gamma_{s'})\qquad\text{if}\ s'=1,\ldots, s-1, s+r+1,\ldots, \ell \\
b\rho(\gamma_{s'})b^{-1}\qquad\text{if}\ s'=s,\ldots, s+r-1\ .
\end{cases}
\end{equation}
Hence, the only conjugacy class which is potentially changed is that of
$\rho(\gamma_{s+r})$.  One can easily show that any deformation
of a Lagrangian representation of the form \eqref{E:bending}
with $b\in O_{L_s}$ is necessarily Lagrangian and coincides with
$\tilde\varphi(\lambda_b)$.

\begin{lem} \label{L:bending}
Let ${\mathcal B}_\rho(s,r)
\subset T_\rho\LHom(\Gamma_\ell, U(n))$ denote the subspace tangent
to the bending deformations of $\rho$ of length $r$ about $L_s$.
Then the $(s+r)$-th component
$
\left[{\mathcal B}_\rho(s,r)\right]_{s+r}
$
is the orthogonal projection of ${\mathfrak o}_s$ to ${\mathfrak
o}_{s+r}^\perp$.
\end{lem}

\Pf
Using \eqref{E:bending} and the calculation in the proof of Lemma \ref{inj}, we see
 that for an
infinitesimal bending $\dot b=B\in{\mathfrak o}_s$,
$$
X_{s'}=\begin{cases}
0\qquad\text{if}\ s'=1,\ldots,s-1,s+r+1,\ldots,\ell \\
({\bf I}-\Ad_{\rho(\gamma_{s'})})B \qquad\text{if}\ s'=s,\ldots, s+r-1 \\
({\bf I}-\Ad_{\sigma_{s+r}})B \qquad\text{if}\ s'=s+r \ .
\end{cases}
$$
Hence, 
$
\left[{\mathcal B}_\rho(s,r)\right]_{s+r}=
\im({\bf I}-\Ad_{\sigma_{s+r}})\bigr|_{\mathfrak o_s}
$,
and the result follows from Lemma \ref{Ad}.
\EPf

Our goal  is to show that twistings and real bendings sweep out the full space of 
 deformations of conjugacy classes
 in a neighborhood of an irreducible Lagrangian representation. 
 We now prove

\begin{prop}[cf.\ Proposition
\ref{maximalrank1}]  \label{maximalrank2}
Let $\rho:\Gamma_\ell\to U(n)$ be an irreducible Lagrangian representation with
$\pi(\rho)={\mathfrak a}\in{\mathcal L}_{I,\ell}(n,{\mathfrak m}, z)$.
 Then near $\rho$, 
$\LRep^{irr.}_{\mathfrak a}(\Gamma_\ell, U(n))$
is a smooth manifold of dimension
$$
\dim\left(\LRep^{irr.}_{\mathfrak a}(\Gamma_\ell, U(n))\right)=
\frac{(\ell-2)}{2}n^2+1-\frac{1}{2}\sum_{s=1}^\ell\sum_{j=1}^{l_s} (\mu^s_j)^2\ .
$$
Moreover, the spectral projection
$$
\pi: \LRep^{irr.}(\Gamma_\ell, U(n))\cap \pi^{-1}\left( 
{\mathcal L}_{I,\ell}(n,{\mathfrak m}, z)\right)
\longrightarrow
{\mathcal P}_{I,\ell}(n,{\mathfrak m}, z)\ ,
$$
is locally surjective and is a fibration near $\rho$.
\end{prop} 

\Pf
As in the proof of Proposition \ref{maximalrank1}, we will first 
concentrate on deformations of
$\rho(\gamma_1)$ up to conjugation.  Thus, we consider bending deformations of $L_s$ and
length $r=\ell-s+1$, for $s=2,\ldots, \ell$.  We also add twist deformations.
It follows from Lemmas \ref{L:twist} and \ref{L:bending} that for
$\rho=\tilde\varphi(\lambda)$,
$$
P_1^\perp{\mathfrak o}_2 +\cdots +
P_1^\perp{\mathfrak o}_\ell
+i\R+\im({\bf I}-\Ad_{\rho(\gamma_1)})\subset [\im D\tilde \varphi_\lambda]_1\ ,
$$
where $P_1^\perp$ is the orthogonal projection to ${\mathfrak
o}_1^\perp$.
Since we are assuming $\rho$ is irreducible, it follows as in the proof of Lemma \ref{inj}
that
$$
({\mathfrak o}_1+\cdots+{\mathfrak o}_\ell)^\perp=\Ker (D\tilde\varphi_\lambda)=i\R\ .
$$
Hence, denoting the traceless part with a subscript $0$,
\begin{align*}
\bigl\{ P_1^\perp{\mathfrak o}_2 +\cdots +
P_1^\perp{\mathfrak o}_\ell
+i\R &+\im({\bf I}-\Ad_{\rho(\gamma_1)})
\bigr\}^\perp 
\\
&=\bigl\{P_1^\perp({\mathfrak o}_1 +\cdots +
{\mathfrak o}_\ell)
+i\R+\im({\bf I}-\Ad_{\rho(\gamma_1)})
\bigr\}^\perp \\
&=\bigl[\{ {\mathfrak o}_1+({\mathfrak o}_1 +\cdots +
{\mathfrak o}_\ell)^\perp\} \cap
{\mathfrak z}(\rho(\gamma_1))
\bigr]_0\\
&=
\left[ {\mathfrak o}_1 \cap {\mathfrak z}(\rho(\gamma_1))
\right]_0
={\mathfrak o}_1\cap {\mathfrak o}_2\ ,
\end{align*}
by Proposition \ref{decomposition} (2).
Since we may do this calculation for any $\rho(\gamma_s)$, and since the variation
preserves the conjugacy classes of $\rho(\gamma_{s'})$, $s'\neq s$ up to the twist
deformations, we have shown that $D\tilde\varphi_\lambda$ is surjective onto
\begin{equation} \label{subspace}
\bigl[ ({\mathfrak o}_1\cap {\mathfrak o}_2)^\perp \times
({\mathfrak o}_2\cap {\mathfrak o}_3)^\perp
\times\cdots\times ({\mathfrak o}_\ell\cap {\mathfrak o}_1)^\perp \bigr]_0
\subset \left[
{\mathfrak z}(\rho(\gamma_1)) \times\cdots\times
{\mathfrak z}(\rho(\gamma_\ell))\right]_0\ ,
\end{equation}
where now the subscript indicates that the sum of the traces vanishes.
By Lemma \ref{L:orth}, this must be exactly the image.
Notice that
$${\mathfrak z}^{ab.}(\rho(\gamma_s))\oplus {\mathfrak
z}^\perp(\rho(\gamma_s))\subset ({\mathfrak o}_s\cap{\mathfrak o}_{s+1})^\perp\ ,$$
 for
all $s$ (cf.\ Proposition \ref{decomposition} and Lemma 
\ref{L:multiplicity}).  Hence, by transversality 
we deduce the local surjectivity and fiber structure onto
the multiplicity space.
We count dimensions:
\begin{itemize}
\item $\ell(n/2)(n+1)$ is dimension of $\Lambda^\ell(n)$;
\item 
By Proposition \ref{decomposition} (4),
$$
\dim{\mathfrak z}(\sigma_s\sigma_{s+1})-(\dim{\mathfrak o}_s\cap {\mathfrak
o}_{s+1})
=(1/2)
\dim{\mathfrak z}(\rho(\gamma_s))+n/2\ .
$$
Hence, the dimension of the subspace in \eqref{subspace} is
$(1/2)\sum_{s=1}^\ell\sum_{j=1}^{l_s} (\mu^s_j)^2 + \ell(n/2)-1$.
\item Finally, $n^2$ is the dimension of $U(n)$-orbit through 
$\rho$ (notice that the action is free;
see also Lemma \ref{inj}).
\end{itemize}
The dimension follows by subtracting the last two items from the first.
This completes the proof.
\EPf

Proposition \ref{maximalrank2} implies that, near irreducible
representations, the allowed holonomies for unitary and Lagrangian
representations coincide.  In particular, a chamber either
 has no Lagrangian representations or is
entirely populated by Lagrangians.

\begin{cor}  \label{isolatedchambers}
Let $\Delta\subset {\mathcal U}_{I,\ell}(n,\mathfrak{m},z)$ be a chamber.  
Then $\Delta\cap{\mathcal L}_{I,\ell}(n,\mathfrak{m},z)\neq\emptyset\iff 
\Delta\subset\Lell(n,\mathfrak{m},z)$.
\end{cor}
\Pf    By Remark \ref{chamberstructure}  (2) and Proposition 
\ref{maximalrank2} it follows that $\Delta\cap{\mathcal L}_{I,\ell}(n,\mathfrak{m},z)$ 
is open.  On the other hand this set is also clearly closed in 
$\Delta$; hence, the result.
\EPf

We also have the

\emph{Proof of Theorem \ref{T:lagrangian}}. 
Assume
 $\LRep_{\mathfrak a}^{irr.}(\Gamma_\ell, U(n))$ is not
empty.
Then by Propositions \ref{maximalrank1} and \ref{maximalrank2}, 
it is a smoothly
embedded half-dimensional submanifold of $\Rep_{\mathfrak
a}^{irr.}(\Gamma_\ell, U(n))$.
By Proposition \ref{P:isotropic}, its tangent space is everywhere
isotropic.  The theorem follows.
\EPf


\subsection{Codimension of the reducibles}

In this section, we 
 use Proposition \ref{maximalrank2}  to 
estimate the size of the set of reducible representations.
Since we will only require the result for $\ell=3$, we restrict to this case.
  We begin with the following simple observation.

\begin{lem} \label{L:mult}
Let $\rho:\Gamma_3\to U(n)$ be irreducible with $\pi(\rho)={\mathfrak a}\in
{\mathcal P}_{I,3}(n,{\mathfrak m},z)$.  Then for at least two values of
$s=1,2,3$,
all multiplicities $\mu_j^s=m^s_j-m^s_{j-1}\leq n/2$.
\end{lem}

\Pf
Suppose not.  Then there are two values of $s$, say $s=1,2$, and $j_1$, $j_2$,
such that $\mu^1_{j_1}>n/2$ and $\mu^2_{j_2}>n/2$.  
If $E_1$ is the $\hat\alpha^1_{j_1}$ eigenspace of $\rho(\gamma_1)$
and $E_2$ is the $\hat\alpha^2_{j_2}$ eigenspace of $\rho(\gamma_2)$, 
then both $\rho(\gamma_1)$ and $\rho(\gamma_2)$, and hence also $\rho(\Gamma_3)$, leave invariant the intersection
$E_1\cap E_2$, which is positive dimensional.  This contradicts the assumption of
 irreducibility.
\EPf

\begin{prop} \label{P:codimension}
Let $\Omega\subset \LRep_{\mathfrak a}(\Gamma_3,U(n))$ 
be an  open connected subset 
containing an irreducible representation. 
Then the set of reducibles $\Omega\cap\LRep_{\mathfrak a}^{red.}(\Gamma_3,U(n))$ 
has codimension $\geq n$.
\end{prop}

\Pf
Suppose ${\mathfrak a}\in {\mathcal P}_{I,3}(n,{\mathfrak m},z)$.
If $\tilde \rho\in \LRep_{\mathfrak a}(\Gamma_3,U(n))$ is reducible,
then we can decompose it into its irreducible components 
$\rho_i$, $i=1, \ldots,
k$, $k\geq 2$. Without loss of generality, we may assume $\rho_i$ and $\rho_j$
are non-isomorphic for $i\neq j$.
  Write: $\pi(\rho_i)={}_i{\mathfrak a}=({}_i\alpha^s_j)\in
{\mathcal P}_{I_i,3}(n_i,{}_i{\mathfrak m},z_i)$.
 Conversely, given a decomposition of ${\mathfrak a}$ into 
${}_1{\mathfrak a},\ldots, {}_k{\mathfrak a}$, 
it suffices to compute the codimension
of the set of all reducibles with $\pi_i(\rho)={}_i{\mathfrak a}$.
We therefore assume this fixed decomposition, and let
 $\cod$ be the codimension of all  reducibles compatible with the decomposition.

For each $s$ let $\mu_j^s$, $j=1,\ldots, l_s$ denote the multiplicities from the
partition $m^s$, and let $\hat\alpha^s_j$ denote the distinct entries of
$\alpha^s$.  We define ${}_i\mu^s_j$ to be the multiplicity of 
$\hat\alpha^s_j$ if it
appears in  ${}_i\alpha^s$, and we set it to zero
otherwise.  The following are easy consequences of this definition.
\begin{align}
\mu^s_j &= \sum_{i=1}^k {}_i\mu^s_j\ , \label{mu}\\
n_i &= \sum_{j=1}^{l_s} {}_i\mu^s_j\ , \label{ni} \\
n&=\sum_{i=1}^k n_i=\sum_{j=1}^{l_s} \mu^s_j\ .\label{n}
\end{align}

Counting dimensions as in the proof of Proposition \ref{maximalrank2} we find
\begin{align}
\cod &=
(3/2)n^2-\bigl((1/2)
\sum_{s=1}^3\sum_{j=1}^{l_s}(\mu^s_j)^2-1\bigr)-n^2
\label{irr}\\
&\qquad-\biggl\{
\sum_{i=1}^k
\bigl[(3/2)n_i^2-((1/2)
\sum_{s=1}^3\sum_{j=1}^{l_s}({}_i\mu^s_j)^2-1
)-n_i^2\bigr]-\bigl( n^2-\sum_{i=1}^k n_i^2\bigr)\biggr\} 
\label{red}\\
&=(1/2)\sum_{s=1}^3
\biggl\{n^2-\sum_{j=1}^{l_s}(\mu^s_j)^2
-\sum_{i=1}^k
\bigl( n_i^2-\sum_{j=1}^{l_s}({}_i\mu^s_j)^2\bigr)\biggr\}\label{cod} +1-k\ . 
\end{align}
The line \eqref{irr} is the dimension count for the irreducibles.  In line
\eqref{red}, we take this dimension for each irreducible factor, and then divide out
by the part of the $U(n)$  which changes the splitting.
It follows that for each $s$ we need to estimate
$$
C_s=
n^2 -\sum_{i=1}^k n_i^2
-\sum_{j=1}^{l_s}(\mu^s_j)^2
+\sum_{i=1}^k\sum_{j=1}^{l_s}({}_i\mu^s_j)^2\ .
$$
Using \eqref{n} we have
$$
n^2=\bigl(\sum_{i=1}^k n_i\bigr)^2=
\sum_{i=1}^k n_i^2+\sum_{i\neq i'}n_in_{i'}\ .
$$
Applying \eqref{ni} to the second term on the right hand side above,
\begin{equation} \label{neq}
n^2-\sum_{i=1}^k n_i^2=
\sum_{i\neq i'}\sum_{j,j'}({}_i\mu^s_j)( {}_{i'}\mu^s_{j'})
=\sum_{i\neq i'}\sum_j({}_i\mu^s_j)( {}_{i'}\mu^s_{j})
+\sum_{i\neq i'}\sum_{j\neq j'}({}_i\mu^s_j)( {}_{i'}\mu^s_{j'})\ .
\end{equation}
On the other hand, from \eqref{mu} we have
\begin{equation} \label{mueq}
\sum_{j=1}^{l_s} (\mu^s_j)^2=
\sum_{i,i'}\sum_{j}({}_i\mu^s_j)( {}_{i'}\mu^s_{j})
=\sum_{i, j}({}_i\mu^s_j )^2
+\sum_{i\neq i'}\sum_{j}({}_i\mu^s_j)( {}_{i'}\mu^s_{j})\ .
\end{equation}
Combining \eqref{neq} and \eqref{mueq}, we find that
$$
C_s=
\sum_{i\neq i'}\sum_{j\neq j'}({}_i\mu^s_j)( {}_{i'}\mu^s_{j'})\ .
$$
We wish to estimate this quantity from below.
Since there are at least two distinct eigenvalues, it follows that:
$C_s\geq 2$.  By Lemma \ref{L:mult}, for at least two values of $s$
we may assume that $\mu^s_j\leq n/2$ for all $j=1,\ldots, l_s$.  
We estimate $C_s$ in this case.

\medskip\noindent Case 1.  Assume that 
for each $i,j$ where ${}_i\mu^s_j\neq 0$ there are 
$i'\neq i$ and $j'\neq j$ such that ${}_{i'}\mu^s_{j'}\neq 0$.  In this case we
have
\begin{equation} \label{E:est1}
C_s \geq 2
\sum_{i,j}({}_i\mu^s_j)\geq 2n\ ,
\end{equation}
by \eqref{ni} and \eqref{n}.

\medskip\noindent Case 2.  If the condition in Case 1 is not satisfied, 
then there
are  $i_0$, $j_0$ such that ${}_{i_0}\mu^s_{j_0}\neq 0$
and for all $i\neq i_0$, $n_i=1$ and ${}_i\mu^s_{j}=1$ if $j=j_0$ and  zero
otherwise. This is true because
 if $n_i\geq 2$, then the $i$-th block  must have at least two distinct
eigenvalues; in particular, one different from ${}_{i_0}\mu^s_{j_0}$. 
 We also have $n_{i_0}-{}_{i_0}\mu^s_{j_0}=n-\mu^s_{j_0}$, and
$n-n_{i_0}=k-1$.  Now
\begin{align*}
\sum_{i\neq i'}\sum_{j\neq j'}({}_i\mu^s_j)( {}_{i'}\mu^s_{j'})
&=2\bigl(\sum_{j\neq j_0}{}_{i_0}\mu^s_j\bigr)(n-n_{i_0})
+\sum_{i_0\neq i\neq i'\neq i_0} ({}_i\mu^s_{j_0})({}_{i'}\mu^s_{j_0})\\
&=2\bigl(\sum_{j}{}_{i_0}\mu^s_j-{}_{i_0}\mu^s_{j_0}\bigr)(n-n_{i_0})
+(1/2)(n-n_{i_0})(n-n_{i_0}-1) \\
&=2(n_{i_0}-{}_{i_0}\mu^s_{j_0})(n-n_{i_0})
+(1/2)(n-n_{i_0})(n-n_{i_0}-1) \\
&=2(n-\mu^s_{j_0})(n-n_{i_0})
+(1/2)(n-n_{i_0})(n-n_{i_0}-1) \ ,
\end{align*}
where in the third line we have used \eqref{n}.   Using
 the assumption that $\mu^s_j\leq n/2$, we have
\begin{equation} \label{E:est2}
C_s\geq n(k-1)+(1/2)(k-1)(k-2)\ .
\end{equation}
Hence, we have bounds on $C_s$ from Cases 1 and 2 at two of the three values of
$s$, and $C_s\geq 2$ at the third value.
Putting \eqref{E:est1} and \eqref{E:est2}
 into the expression \eqref{cod} we find  three possibilities:
$$
\cod\geq
\begin{cases}
2n+2-k\ ; \\
n+(1/2)\left\{ n(k-1)+(1/2)(k-1)(k-2)\right\}+2-k\ ;\\
n(k-1)+(1/2)(k-1)(k-2)+2-k\ .
\end{cases}
$$
It is easily verified that the quantities on the right are all $\geq n$, with
equality in the last case at $k=2$.
Since this is true for all of the finitely many possible types of reduction, the
proof is complete.
\EPf


\section{Proof of the Main Theorem} \label{S:proof}

We have shown in Proposition \ref{maximalrank2} that 
$\LRep_{\mathfrak a}^{irr.}(\Gamma_\ell, U(n))$, 
if not empty, is a smoothly embedded submanifold of 
$\Rep_{\mathfrak a}^{irr.}(\Gamma_\ell, U(n))$.  In this section, we prove
the existence of a Lagrangian representation with given holonomy
whenever a unitary representation with the same holonomy exists. 
  We first reduce the problem to the case of triples.
  
  \begin{prop}  \label{triple}
  Suppose Theorem \ref{maintheorem} holds for $\ell=3$.  Then it holds for all $\ell$.
  \end{prop}
  
  \Pf  By induction.  Assume  Theorem \ref{maintheorem} holds for some $\ell\geq 3$, and
also for $\ell=3$. 
 We show that it also holds for $\ell+1$. 
 Let $A_1,\ldots, A_{\ell+1}$ be
 unitary matrices satisfying $A_1\cdots A_{\ell+1}={\bf I}$ with 
given spectra.  
 By induction, we may find Lagrangians $L_1,\ldots, L_{\ell-1}$ such that
$\spec(A_i)=\spec(\sigma_{L_{i-1}}\sigma_{L_i})$, $i=1,\ldots, \ell-1$, and
$\spec(A_\ell A_{\ell+1})=\spec(\sigma_{L_{\ell-1}}\sigma_{L_0})$, where $L_0$ is as in Section \ref{S:linalg}.
Write:  $B_1 B_2 B_3={\bf I}$,
where $B_1\sim A_{\ell+1}^{-1}$, 
$B_2\sim A_{\ell}^{-1}$,  and $B_3=\sigma_{L_{\ell-1}}\sigma_{L_0}$.
Using the result for $\ell=3$ we may find Lagrangians $L'$, $L''$ such that $B_1\sim
\sigma_{L_0}\sigma_{L'}$, 
$B_2\sim \sigma_{L'}\sigma_{L''}$, and $B_3\sim \sigma_{L''}\sigma_{L_0}$. By
Lemma \ref{L:sym}, 
both $\sigma_{L_{\ell-1}}\sigma_{L_0}$ and $\sigma_{L''}\sigma_{L_0}$ are
conjugate by elements in $O(n)$ to  diagonal matrices. 
 Since they furthermore have the same spectrum,
it follows from Proposition \ref{P:tau1} that there 
is some $g\in O(n)$ with $gL''=L_{\ell-1}$.
Set $L_\ell=gL'$.  Then $A_\ell\sim \sigma_{L_{\ell-1}}\sigma_{L_\ell}$, and 
 $A_{\ell+1}\sim \sigma_{L_{\ell}}\sigma_{L_0}$,  and the result follows.
  \EPf
  
  By Proposition \ref{triple}, it suffices to prove
 Theorem \ref{maintheorem} for 
triples of Lagrangians.  For the rest of this section, we consider the 
problem of specifying three conjugacy classes.   To simplify notation, 
we will omit the subscript ``$\ell=3$", and write $\Gamma$ for $\Gamma_3$, 
and $\overline\U^\ast(n)$ for $\overline{\mathcal U}^\ast_{I,3}(n)$, for example.

  \begin{dfn} \label{relred}
  A reducible representation 
$\rho : \Gamma \to U(n_1)\times\cdots\times U(n_k)\hookrightarrow 
U(n)$, $\sum_{i=1}^k n_i=n$,  will
be called \emph{relatively irreducible} with respect to 
$U(n_1)\times \cdots\times U(n_k)$ if the induced representations 
$\rho_i: \Gamma \to U(n_i)$ are irreducible for each $i=1,\ldots, k$.
  \end{dfn}

Our goal is to show that $\overline{\mathcal L}^\ast_I(n)=\overline{\mathcal
U}^\ast_I(n)$, for all $I$ and $n$.
Using the stratification of $\overline{\mathcal P}_I^\ast(n)$
described in Section \ref{S:conjugacy},
  the argument proceeds by induction on the four parameters available:
  \begin{itemize}
  \item
  Fix the rank $n$.   We assume that we have shown 
$\L(\tilde n,{\mathfrak m},z)=\U(\tilde n,{\mathfrak m}, z)$ 
for all $\tilde n<n$ and all $({\mathfrak m},z)$.  
The result for $U(1)$ or $U(2)$ representations holds, 
as has already been mentioned.
  \item Next, fix a multiplicity structure ${\mathfrak m}$.
 Assume we have proven that 
  $\L(\tilde n,{\mathfrak p},z)=\U(\tilde n,{\mathfrak p}, z)$ 
for all ${\mathfrak p}<{\mathfrak m}$ and all $z$.  We may clearly 
do this, since a partition giving multiplicity $n$ for each $s$ 
corresponds to $U(1)$ representations. 
  \item   Fix a subset $z\subset\{1,2,3\}$ and assume that 
$\L( n,{\mathfrak m},\tilde z)=\U(n,{\mathfrak m}, \tilde z)$ for all 
$z\varsubsetneq\tilde z$.  We will justify this assumption below.  
  \item Finally, the last part of the inductive scheme is to assume 
that ${\mathcal L}_{\overline I}( n, \bar{\mathfrak m},\bar z)=
\U(n,\bar{\mathfrak m},\bar z)$
 for all $\overline I<I$, and all $\bar{\mathfrak m}$ and $\bar z$.  
Notice that $\overline I=0$ involves only the trivial representation.
  \end{itemize}
  
  \noindent
If the stratum $\P(n,{\mathfrak m},z)$ is degenerate, 
then either $\U(n,{\mathfrak m},z)=\emptyset$, in which case there is 
nothing to prove, or each $\rho$ with 
$\pi(\rho)\in\U(n,{\mathfrak m},z)$ is reducible by Proposition 
\ref{maximalrank1}.  
Hence, by induction on the rank $n$, $\L( n,{\mathfrak m},z)=
\U( n,{\mathfrak m}, z)$ if $\P(n,{\mathfrak m},z)$ is degenerate.
Thus, we assume that $\P(n,{\mathfrak m},z)$ is nondegenerate.  
If $\L(n,{\mathfrak m},z)\neq\U(n,{\mathfrak m},z)$ then there is a connected component $\Delta$ of 
  $\U(n,{\mathfrak m},z)\setminus\L(n,{\mathfrak m},z)$ 
which by Corollary \ref{isolatedchambers} is a union of chambers.  
By Remark \ref{chamberstructure} (1), $\partial\Delta$ consists of a 
union of 
convex subsets of affine planes.  By Proposition \ref{maximalrank2}, 
it follows that any $\rho\in\LHom(\Gamma, U(n))$ for 
which $\pi(\rho)\in\partial \Delta$ is reducible. 
 Finally,  we claim that $\partial\Delta\cap\intU(n,{\mathfrak m},z)$
 is unbounded. 
 To see this, choose
$\rho\in\partial \Delta\cap\intU(n,{\mathfrak m},z)$ contained in a cell
of minimal dimension.  Then $\rho$ is relatively irreducible with
respect to some reduction $U(n_1)\times\cdots\times U(n_k)$ (see
Definition \ref{relred}).  Among the
induced representations $\Gamma\to U(n_j)$ there must be one, say
$\rho_j$, that is
nontrivial, since the total index is positive.
Hence, 
$\pi(\rho_j)\in
{\overset{\circ}{\mathcal U}}_{I_j}(n_j, {}_j{\mathfrak m}, z_j)$ for
some induced multiplicity stucture.  Since
${\overset{\circ}{\mathcal U}}_{I_j}(n_j, {}_j{\mathfrak m}, z_j)$ is
positive dimensional,  the claim follows from this fact.

 From the discussion above and the description
 of the stratification and wall structure
 in Sections \ref{S:conjugacy} and \ref{S:inequalities} we see
 that there are 
four (not necessarily exclusive) possibilities:
  \begin{enumerate}
  \item $\partial\Delta$ intersects an outer wall in 
$\P(n,{\mathfrak m},z)$;
  \item $\partial\Delta$ intersects a stratum 
$\P(n,{\mathfrak p},z)$, ${\mathfrak p}<{\mathfrak m}$;
  \item $\partial\Delta$ intersects a stratum $\P(n,{\mathfrak p},
\tilde z)$, ${\mathfrak p}\leq {\mathfrak m}$, $z\varsubsetneq\tilde z$;
  \item $\partial\Delta$ intersects a stratum 
${\mathcal P}_{\overline I}(n,\bar{\mathfrak p},\bar z)$, 
for some $\overline I < I$, $z\subset \bar z$.
  \end{enumerate}
  
  \noindent
  In each case, our inductive hypothesis assumes the result for the
 lower dimensional stratum, and we will use this below to derive a 
contradiction.  Here we remark that possibility (3) does not occur
 if $z=\{1,2,3\}$.  The derivation of a contradiction for 
this case therefore justifies the inductive hypothesis on $z$.
  The structure of the  argument deriving a contradiction is actually 
identical for each of the four possibilities above, 
\emph{mutatis mutandis}.  We will give a detailed account of how 
this works in case (1),  the modifications necessary for
 the other cases being straightforward.

Consider then  the case where 
$\partial {\Delta}$ intersects the outer walls 
$\W(n,{\mathfrak m},z)$ at a point in $\P(n,{\mathfrak m},z)$.  
To simplify notation, for the following discussion we set
$\U=\U(n,{\mathfrak m},z)$, $\W=\W(n,{\mathfrak m},z)$, 
$\P=\P(n,{\mathfrak m},z)$, and $\Lambda_I=
\Lambda^3_I(n,{\mathfrak m},z)$.
Also, let $l_s$ be the lengths of the partitions $m^s$, $s=1,2,3$.
The intersection $H=\partial \Delta\cap \W$ is a union of convex 
subsets of intersections of affine planes corresponding to 
reductions of Lagrangian representations.
We claim that $H$ must have positive codimension in $\W$.  For if not, we could find
an outer wall $W$ and  point ${\mathfrak a}\in \partial\Delta\cap W$ such that
${\mathfrak a}\not\in W'$ for any outer \emph{or inner} wall $W'$.  In particular, if
$N$ is a sufficiently small neighborhood of $\mathfrak a$, then $N\cap \intU
=N\cap \Delta$.    By the induction hypothesis, we may find (a reducible)
$\rho\in {\mathcal L}_I$ such that $\pi(\rho)={\mathfrak a}$.   Now any Lagrangian
may be perturbed slightly to give an irreducible Lagrangian representation $\tilde
\rho$.  It follows from Proposition \ref{P:structure} (1) that for sufficiently small
perturbations, $\pi(\tilde \rho)\in N\cap\intU\subset\Delta$;
contradiction.

Hence, we may assume $H$ has positive codimension.
To illustrate the basic idea of the proof, suppose first that $H$ has 
codimension one inside $\W$, so that $H$ locally disconnects $\W$.
We  choose ${\mathfrak a}\in H$  with minimal valency with respect 
to the outer wall structure along $H$.  By this we mean that there 
are outer walls $W_1,\ldots, W_p$ meeting at ${\mathfrak a}$, 
and $p\geq 1$ is the minimal number of such intersections among all 
points in $H$.
 With this choice,  and using the 
convexity of $\U$, we see that  the number $p$ of outer walls meeting at 
${\mathfrak a}$ is $1$ or $2$.  Let us first assume that $p=1$, and 
let $W$ denote the outer wall in question.  Choose a neighborhood $U$ 
of ${\mathfrak a}$ in the wall $W$ such that $H\cap U$ is a cell.  
Since $W$ is the only outer wall at ${\mathfrak a}$, we may also 
assume that the neighborhood $U$ is contained in $\U$. 
Let $N$ be a neighborhood of ${\mathfrak a}$ in $\P$ such that the 
following hold:
\begin{enumerate}
\item $U=N\cap W$;
\item  $N\setminus W$ consists precisely of two components $N^+$, $N^-$;
\item $N^-\cap \U=\emptyset$ and $N^+\subset \U$ is homeomorphic to a ball;
\item $N^+\setminus \Delta$ has the topology of $U\setminus H$.
\end{enumerate}
 Choose a point $\rho\in \pi^{-1}(W)$ as follows:  $W$ corresponds to a
 reduction $U(k)\times U(n-k)$. 
 We may find a point $\rho$, $\pi(\rho)=
{\mathfrak a}$, such that $\rho$ is relatively irreducible with respect to
$U(k)\times U(n-k)$. 
By Proposition \ref{maximalrank2}, we may assume that $\Lambda_I$ is a
manifold near $\rho$.   With this understood, 
let $\widetilde B\subset \Lambda_I$  be a ball about $\rho$ such that 
$\pi(\widetilde B)\cap W\subset U$.  By our choice of $\rho$ it follows, again by
Proposition \ref{maximalrank2},  that 
$\pi(\widetilde B) $ intersects \emph{both} components of 
$U\setminus H$.  
By Proposition \ref{P:codimension}, $\widetilde B\cap\Lambda^{irr.}_I$ 
is connected; hence, so is $\pi(\widetilde B\cap\Lambda^{irr.}_I)$.  
On the other hand,  by the previous remark, 
$\pi(\widetilde B\cap\Lambda^{irr.}_I)\subset N^+\setminus \Delta$ 
must intersect both components of $N^+\setminus \Delta$.  
This contradicts the connectedness of 
$\pi(\widetilde B\cap\Lambda^{irr.}_I)$ (see Figure 1).

\setlength{\unitlength}{1cm}
\begin{picture}(14,8)

\put(4.75,4.75){$W$}
\put(7.5,5.3){$H$}
\put(6.5,2){$\Delta$}
\ifx\pdftexversion\undefined
\put(3,1){
{\scalebox{.6}{\includegraphics{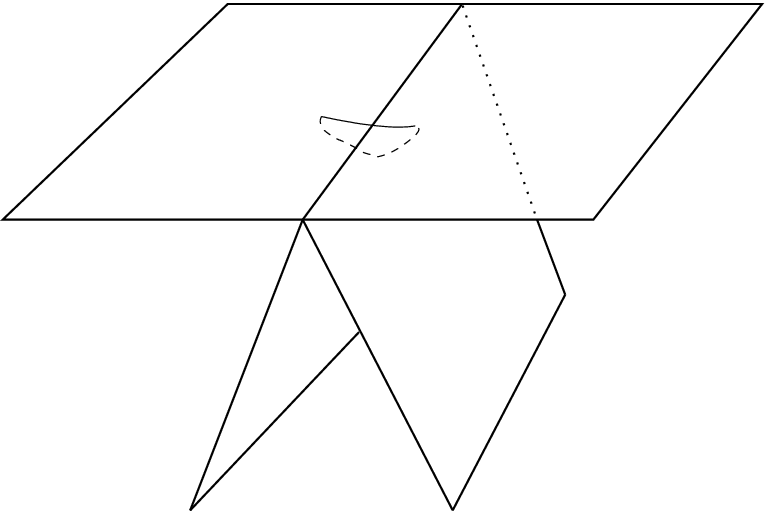}}}}
\else
\put(3,1){
{\scalebox{.6}{\includegraphics{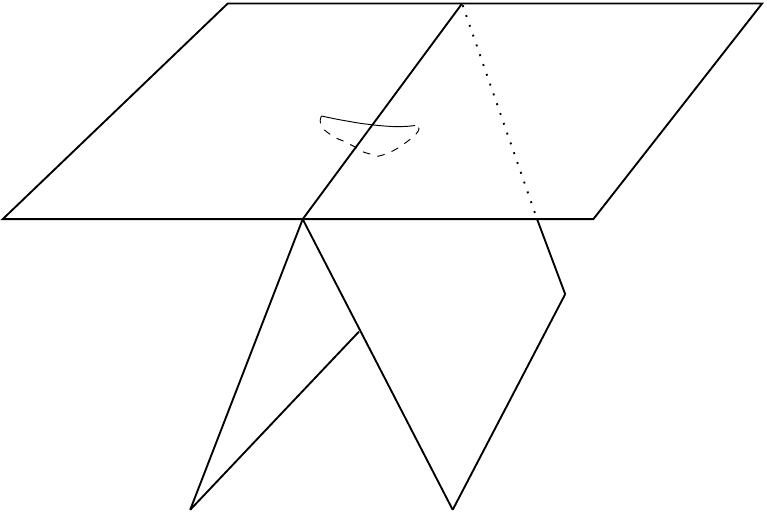}}}}
\fi
\put(2,0){Figure 1: Intersection $H$ of a chamber $\Delta$ with an outer wall $W$}
\end{picture}

\medskip

The case $p=2$ requires only a small modification of the above 
argument:  Let $W_1$ and $W_2$ be outer walls meeting along $H$ at 
${\mathfrak a}$.  We choose the set $U\subset W_1\cup W_2$ to consist 
of two pieces: $U_1=U\cap W_1\subset W_1\cap\U(n)$, and
$U_2=U\cap W_2\subset W_2\cap\U(n)$.  Since ${\mathfrak a}$ is at the 
intersection of precisely two outer walls, it corresponds to a 
reduction of the form $U(n_1)\times U(n_2)\times U(n_3)$;  
the wall $W_1$ corresponds to a $U(n_1+n_2)\times U(n_3)$ reduction, 
say, and 
the wall $W_1$ corresponds to a $U(n_1)\times U(n_2+n_3)$ reduction.  
Now since deformations along the wall $W_1$ can only take values on 
one side of $W_2$, and vice-versa, it follows that the image by $\pi$ 
of a neighborhood of any $\rho$, $\pi(\rho)={\mathfrak a}$,  
intersects \emph{both} components of $U\setminus H$. In the choice of 
the neighborhood $N$ we modify the first two criteria so that
\begin{enumerate}
\item[1'.] $U=N\cap (W_1\cup W_2)\cap \U(n)$;
\item[2'.]  $N\setminus (W_1\cup W_2)\cap \U(n)$ consists precisely of two components $N^+$, $N^-$,
\end{enumerate}
and keep items (3) and (4) as above.  
  The rest of the argument then proceeds exactly as before.

Next, let us consider the case where $H$ has higher codimension $d$, 
$d\geq 2$, in $\W(n)$.  
If we again choose ${\mathfrak a}\in H$  with minimal valency with 
respect to the outer wall structure along $H$, then we see that at most $d+1$ outer walls meet at ${\mathfrak a}$.
As before, we first consider the case where there is just one
 outer wall $W$.  Choose a neighborhood $U$ of 
${\mathfrak a}$ in $W$  as above.  We also choose $N$ satisfying conditions (1-4) above.
Let $D\subset U$ be a cell in $U$ of dimension equal to the
 codimension $d$ of $H$ in $W$ and intersecting $H$ precisely
 in $\mathfrak a$.  Hence, the boundary $\partial D$ is the link of $H$ in
 $W$.   We regard $D$ as the image of a continuous map, 
$f: B^d\longrightarrow U$.  We may further assume that
 $f=\pi\circ\tilde f$ for a map $\tilde f:B^d\longrightarrow
 \Lambda_I$, taking the origin to $\rho$.
Indeed, choosing a relatively irreducible $\rho$ and using  Proposition \ref{maximalrank2},
$\pi: \pi^{-1}(W)\cap\Lambda_I\to W$
 is a fibration in a neighborhood of $\rho$ and $\pi(\rho)=\mathfrak a$.  
Hence, we may define $\tilde f$ by taking a section of this fibration.
\medskip

\noindent {\bf Claim:}\ $\dim H\geq \sum l_s - n-|z|$.  

\Pf
Assume that $\rho$ is relatively irreducible with respect to a
reduction $U(n_1)\times\cdots\times U(n_k)$.
 Then restricted to
representations near $\rho$ which are relatively irreducible of this
type, the map $\pi$ is locally surjective onto 
$ \overline{\mathcal P}_{I_1}(n_1)\times  \cdots\times 
\overline{\mathcal P}_{I_k}(n_k)$ (cf.\ Proposition \ref{maximalrank1}).
 Assume first that $|z|\neq 3$.  Then all $I_j>0$. In particular,
$$
\dim\left(\overline{\mathcal P}_{I_1}(n_1)\times  \cdots\times 
\overline{\mathcal P}_{I_k}(n_k)\right)=
\sum_{j=1}^k (3 n_j -1) = 3 n-k\ .
$$
Since $l_s$ is the number of distinct eigenvalues of
$\rho(\gamma_s)$,  
it follows that
$$
\dim H = 3 n-k -\sum_{s=1}^3  (n-l_s)-|z|
= \sum_{s=1}^3 l_s-k-|z|\geq 
 \sum_{s=1}^3 l_s-n-|z|\ .
$$
Now suppose that $I_1=\cdots = I_q=0$ for some $1\leq q<k$,  and
 $I_j\neq 0$ for $j=q+1,\ldots, k$.   Since we are assuming
$\pi(\rho)\in\P(n,{\mathfrak m},z)$, 
this can only happen if $z=\{1,2,3\}$, i.e.\ $|z|=3$.
Also, $n_1=\cdots=n_q=1$.  It follows that
\begin{align*}
\dim\left(\overline{\mathcal P}_{I_1}(n_1)\times  \cdots\times 
\overline{\mathcal P}_{I_k}(n_k)\right)&=
\dim\left(\overline{\mathcal P}_{I_{q+1}}(n_{q+1})
\times  \cdots\times 
\overline{\mathcal P}_{I_k}(n_k)\right) \\
&= \sum_{j=1}^{k-q} (3 n_{q+j}-1)=3(n-q)-(k-q)\ .
\end{align*}
Now for each $s=1,2,3$, either $q=m^s_1$, in which case there are
precisely $l_s-1$ distinct nonzero eigenvalues among the remaining
$n-q$; or, $q< m^s_1$, in which case there are $l_s$ distinct
eigenvalues, but one of them is zero.  In both cases, this imposes:
$n-q-(l_s-1)$ conditions on the eigenvalues.
Hence, we have
$$
\dim H= 3(n-q)-(k-q) -\sum_{s=1}^3(n-q-(l_s-1)) =  
\sum_{s=1}^3 l_s-(k-q)-3\ .
$$
Since $k-q\leq n-1$, and $|z|=3$, 
 the claim follows in this case as well.  \EPf

Now $d=\dim W  -\dim H\leq \sum_{s=1}^3 l_s -2-|z| -(\textstyle\sum_{s=1}^3
 l_s -n-|z|)= n-2$.
Notice that this computation is still valid even if $\sum_{s=1}^3l_s-n-|z|\leq 0$.
  By Proposition \ref{P:codimension}, $\Lambda^{red.}_I$ has 
codimension at least: $n>n-2\, $ in $\Lambda_I$.
Hence,  we may find a perturbed map $\tilde f_\varepsilon :
 B^d\to \Lambda^{irr.}_I$.  For sufficiently small perturbations
 we clearly may assume that $f_\varepsilon=\pi\circ
 \tilde f_\varepsilon$ has image in $N$.  It follows that in fact 
$f_\varepsilon: B^d\to N^+\setminus \Delta$.  Now $N^+\setminus \Delta$
 has the topology of $U\setminus H$, and under this equivalence 
$f_\varepsilon(\partial B^d)$ is the link of $N\cap\Delta$. 
 The continuous extension of $f_\varepsilon$ to $B^d$ is therefore
 a contradiction.

When the number $p$ of outer walls meeting at ${\mathfrak a}$ is 
greater than one, the configuration of outer walls at 
${\mathfrak a}$ forms a ``corner" in $\W$ (see Figure 2).  As in the case $p=2$ 
above, we want to choose the set $U$ to mimic this configuration.
The technical result we will require is the following:

\begin{lem} \label{walls}
Suppose that $\rho\in \Lambda_I$ is such that $\pi(\rho)$ lies in the
 intersection $W_1\cap\cdots\cap W_p$ of  $p$ distinct outer walls,
 where $p$ is the minimal such number,  and that $\rho$ is relatively 
irreducible with respect to the reduction corresponding to
 $W_1\cap\cdots\cap W_p$.  Then for any small neighborhood
 $\Omega\subset\Lambda_I$ of $\rho$ there is a continuous
 map $\tilde f: B^{p-1}\to\Omega$ satisfying the following:
\begin{enumerate}
\item $\tilde f(0)=\rho$;
\item  $\pi\circ\tilde f(\partial B^{p-1})\subset W_1\cup\cdots 
\cup W_p$;
\item $\pi\circ\tilde f(\partial B^{p-1})\cap W_1\cap\cdots 
\cap W_p=\emptyset$.
\end{enumerate}
Moreover, $\tilde f$ may be chosen to vary continuously with
 $\rho$ satisfying the hypothesis.
\end{lem}

\begin{picture}(14,8)

\put(5.25,4.25){$W_1$}
\put(9.75,3.75){$W_2$}
\put(7.75,2.25){$W_3$}
\put(6.7,3.2){$H$}
\put(9.1,5.1){$\Delta$}
\ifx\pdftexversion\undefined
\put(4,2){
{\scalebox{.6}{\includegraphics{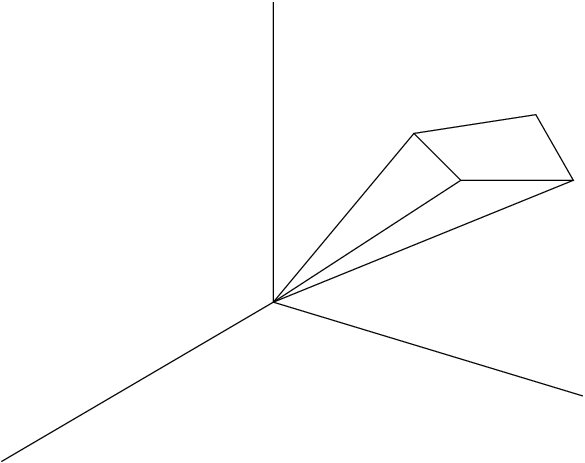}}}}
\else
\put(4,2){
{\scalebox{.6}{\includegraphics{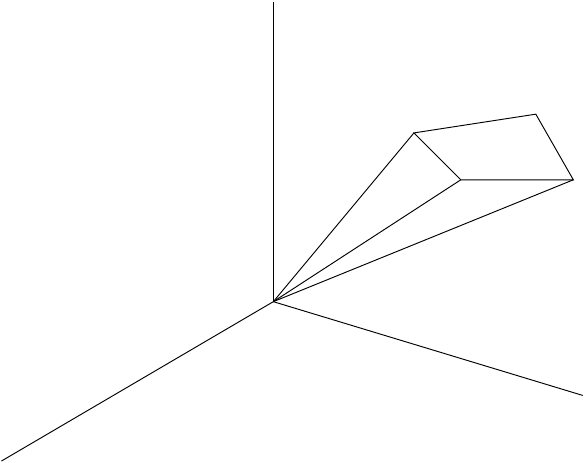}}}}
\fi
\put(2,1){Figure 2: Intersection $H$ of a chamber $\Delta$ with three outer walls}
\end{picture}
\medskip

\noindent
Given the lemma, the rest of the
 argument proceeds as in the previous paragraph.  
Indeed, choose ${\mathfrak a}\in H$  with minimal valency with
 respect to the outer wall structure along $H\subset W_1\cap\cdots
 \cap W_p$, and choose a neighborhood $N$ of ${\mathfrak a}$
 such that $N\setminus (W_1\cup\cdots\cup W_p)\cap \U$ consists
 precisely of two components $N^+$, $N^-$, and which also satisfies 
items (3) and (4) above. Let $\rho$, $\pi(\rho)={\mathfrak a}$, 
be relatively irreducible, and choose a neighborhood $\Omega$ of 
$\rho$ such that $\pi(\Omega)\subset N$.  Choose a continuous map 
$\tilde g: B^{d+1-p}\to \Omega$ such that $\pi\circ \tilde  
g: B^{d+1-p}\to W_1\cap\cdots\cap W_p$ is transverse to $H$ at 
${\mathfrak a}$.  As before, we can do this
 because $\rho$ is relatively irreducible.  Now use  Lemma 
\ref{walls} to extend $\tilde g$ to a continuous map: $\tilde f: 
B^d\simeq B^{d+1-p}\times B^{p-1}\to \Omega$.  By the construction, 
we can easily arrange that
 \begin{equation} \label{boundary1}
 f=\pi\circ\tilde  f(\partial B^{d+1-p}\times\{y\})\cap H=\emptyset\ ,
 \end{equation}
  for all $y\in B^{p-1}$.  By  Lemma \ref{walls} (3) we also have
  \begin{equation} \label{boundary2}
  f(\{x\}\times\partial B^{p-1})\cap W_1\cap\cdots \cap W_p=\emptyset\ ,
  \end{equation}
   for all $x\in B^{d+1-p}$.  It follows from \eqref{boundary1} 
and \eqref{boundary2} that $ f: S^{d-1}\to W_1\cup\cdots \cup W_p$ 
is a link of $H$ in $W_1\cup\cdots \cup W_p$.  We may now perturb the 
map $\tilde f$ as above so that $f_\varepsilon(S^{d-1})\subset 
N^+\setminus \Delta$ is a link of $N^+\cap\Delta$.  
The extension  $f_\varepsilon(B^d)\subset N^+\setminus \Delta$ 
gives a contradiction as before.

 \medskip
 \emph{Proof of Lemma \ref{walls}}.  
Suppose $\rho$ is of type $U(n_1)\times\cdots\times U(n_p)\times 
U(n_{p+1})$, where each wall $W_i$ corresponds to a reduction 
$U(n)\to U(n_i)\times U(n-n_i)$, $i=1,\ldots, p$.  Let 
$\rho=(\rho_1, \ldots, \rho_p,\rho_{p+1})$ be the irreducible factors.
Notice that the  assumption of minimal valency of ${\mathfrak a}$ 
implies that $n_{p+1}=n-\sum_{i=1}^p n_i\neq 0$.  
Let $e_1\cdots e_p$ be a $(p-1)$-simplex in $\R^{p-1}$ with the 
origin $e_0$ as barycenter.  For each $i$, we may find a path 
$\tilde g_i^\prime(t)$ of Lagrangian representations into 
$ U(n_{p+1}+n_i)$ such that $\tilde g_i^\prime(0)=(\rho_i,\rho_0)$
 and $\tilde g_i^\prime(t)$ is irreducible for $t\neq 0$. Keeping 
the other factors fixed, these define paths
$$
\tilde g_i : [0,1]\longrightarrow \LHom(\Gamma,U(n_1)\times\cdots
\times\widehat U(n_i)\times\cdots \times U(n_p)\times 
U(n_{p+1}+n_i))\ ,
$$
where  $\ \widehat{} \ $ means that factor is deleted.  
Combining these paths defines a continuous map
$\tilde f: \cup_{i=1}^p e_0e_i\to\Omega$.  
Suppose inductively that we have defined $\tilde f$ on all 
simplices of the form $e_{i_1}\cdots e_{i_k}$, $2\leq k< p-1$, 
$1\leq i_1<\cdots <i_k$.  For each such simplex, let
  $\{ j_1, \ldots, j_{p-k}\}$ be the complimentary set to 
$\{i_1, \ldots, i_{k}\}$ in $\{1,\ldots, p\}$. 
 We will assume $\tilde f$ has been defined
such that the following hold:
\begin{enumerate}
\item $\pi\circ\tilde f( e_{i_1}\cdots e_{i_k})\subset W_{j_1}\cap\cdots\cap W_{j_{p-k}}$;
\item For each $x\in e_{i_1}\cdots e_{i_k}$, $\tilde f(x)$ is relatively irreducible with respect to the decomposition $U(n_{j_1})\times U(n_{j_{p-k}})\times U\bigl(n-\sum_{\mu=1}^k n_{j_\mu} \bigr)$;
\item $\pi\circ\tilde f( e_{i_1}\cdots e_{i_k})\cap W_1\cap\cdots\cap W_p=\emptyset$.
\end{enumerate}
We now extend $\tilde f$ to a simplex of the 
form $e_{i_1}\cdots e_{i_{k+1}}$ as follows.
By assumption (1), for the complimentary set of indices 
$\{ j_1, \ldots, j_{p-k-1}\}$ we have 
$\pi\circ f(\partial(e_{i_1}\cdots e_{i_{k+1}}))\subset 
W_{j_1}\cap\cdots\cap W_{j_{p-k-1}}$.  Assuming $\Omega$ has been 
chosen sufficiently small so that $\Omega\cap 
\pi^{-1}(W_{j_1}\cap\cdots\cap W_{j_{p-k-1}})$ is contractible, 
we may extend $\tilde f$ to a map $e_{i_1}\cdots e_{i_{k+1}}\to 
W_{j_1}\cap\cdots\cap W_{j_{p-k-1}}$.  Applying the same codimension 
argument we have used several times already, we can further assume 
that this extended map satisfies conditions (2) and (3) as well.  
Continuing in this way, we have defined $\tilde f$ on the boundary of 
$e_1\cdots e_p$. 
  Recall that  $f$ is also defined on the one simplices $e_0e_i$, $i=1,\ldots, p$.
Again using contractibility of $\Omega$, we extend $\tilde f$ 
inductively and arbitrarily to simplices of the form $e_0e_{i_1}\cdots
 e_{i_k}$, $k=1,\ldots, p$.  
This completes the definition of $\tilde f$.
\EPf


\section{Examples} \label{S:examples}

In this last section, we illustrate some of the ideas in the paper by explicity giving
the wall structure for the cases: $\ell=3$, $n=2,3$. For convenience, 
we will only consider distinct eigenvalues different from unity. 
The case of $U(2)$ representations
was first proven \cite{JW}, and more generally \cite{Bi1}. 
 The inequalities were later derived from
spherical triangles in \cite{FMS}.

Let us first introduce some useful notation. For integers $i_s$, $1\leq i_s\leq n$,
$s=1,\ldots, \ell$, define the collection of subsets as in Section \ref{S:inequalities}
$ \wp_{(1)}=(\wp_{(1)}^s)$, $\wp_{(1)}^s=\{i_s\}$.
For ${\mathfrak a}=(\alpha_j^s)\in \overline{\mathcal A}_\ell(n)$, we will use the
notation (cf.\ \eqref{relativeindex})
$$
[i_1,\ldots, i_s]_{\mathfrak a}=I({\mathfrak a},
\wp_{(1)})=\sum_{s=1}^\ell\alpha^s_{i_s}\ .
$$ 
By a permutation of $[i_1,\ldots, i_s]_{\mathfrak a}$, 
we mean a quantity of the form:
$[i_{\tau(1)},\ldots, i_{\tau(s)}]_{\mathfrak a}$,
for some $\tau$ in the  group of permutations of  $\{1,\ldots,\ell\}$.
With this understood, we may write the $U(2)$ inequalities as

\begin{thm}[cf.\ \cite{Bi1}, \cite{FMS}]
There exist representations $\rho: \Gamma_3\to U(2)$ with 
${\mathfrak a}=\pi(\rho)\in {\mathcal U}_{I,3}(2)$, if and only if
\begin{itemize}
\item  $I=2$, and:  $[2,1,1]_{\mathfrak a}\leq 1$,  plus all permutations;
\item  $I=3$, and: $[2,2,1]_{\mathfrak a}\leq 2\leq [2,2,2]_{\mathfrak a}$, 
plus all permutations; or,
\item  $I=4$, and: $[2,1,1]_{\mathfrak a}\leq 2$, plus all permutations;
\end{itemize}
\end{thm}

\noindent The bounds on the index come from Proposition \ref{indexbounds1}.
Notice that for each index there are no inner walls.
Indeed,   any equality of the form: $[i_1,i_2,i_3]_{\mathfrak a}=K$
implies
$$
I=[i_1,i_2,i_3]+[\bar i_1,\bar i_2,\bar i_3]=K+ [\bar i_1,\bar i_2,\bar i_3]
\geq K+1\ ,
$$
where $\bar i_s=\{1,2\}\setminus\{i_s\}$.  Now if $I=2$,  for example, then $K=1$,
and it is easy to see that the outer walls are the only possible solutions for distinct
nonzero eigenvalues.

We have used a duality in the wall structure.  In general, if $\wp_{(k)}=(\wp_{(k)}^s)$
is a collection of subsets of $\{1,\ldots, n\}$ of cardinality $k$, then let
$\wp_{(k)}^\ast$ denote the collection of subsets of cardinality $n-k$ defined by
$(\wp_{(k)}^\ast)^s=(\wp_{(k)}^s)^c$.  It follows that
$
I({\mathfrak a}, \wp_{(k)})+ I({\mathfrak a}, \wp_{(k)}^\ast)
=I({\mathfrak a})
$.
So an inequality of the form
$I({\mathfrak a}, \wp_{(k)})\leq K$ may be written
$ I({\mathfrak a}, \wp_{(k)}^\ast)\geq 
I({\mathfrak a}) -K$.  In particular, this means that for $n=3$
we may express all the inequalities in terms of the $[i_1,\ldots, i_\ell]_{\mathfrak
a}$'s.

\begin{thm}
There exist representations $\rho: \Gamma_3\to U(3)$ with 
${\mathfrak a}=\pi(\rho)\in {\mathcal U}_{I,3}(3)$, if and only if
\begin{itemize}
\item  $I=3$,  and
\begin{align*}
[3,1,1]_{\mathfrak a}\, ,\ [2,2,1]_{\mathfrak a}\leq &1\leq
[3,3,1]_{\mathfrak a}\, ,\ [3,2,2]_{\mathfrak a}\\
&2\leq [3,3,2]_{\mathfrak a}
\end{align*}
plus all permutations;
\item  $I=4$,  and
\begin{align*}
[2,1,1]_{\mathfrak a}\leq &1\leq
[3,2,1]_{\mathfrak a}\, ,\ [2,2,2]_{\mathfrak a}\\
[3,3,1]_{\mathfrak a}\, ,\ [3,2,2]_{\mathfrak a}\leq
&2\leq [3,3,3]_{\mathfrak a}
\end{align*}
plus all permutations;
\item  $I=5$,  and
\begin{align*}
[1,1,1]_{\mathfrak a}\leq &1\leq
[2,2,1]_{\mathfrak a}\, ,\ [3,1,1]_{\mathfrak a}\\
[3,2,1]_{\mathfrak a}\, ,\ [2,2,2]_{\mathfrak a}\leq
&2\leq [3,3,2]_{\mathfrak a}
\end{align*}
plus all permutations; or,
\item  $I=6$,  and
\begin{align*}
 &1\leq
[2,1,1]_{\mathfrak a}\\
[3,1,1]_{\mathfrak a}\, ,\ [2,2,1]_{\mathfrak a}\leq
&2\leq [3,3,1]_{\mathfrak a}\, ,\ [3,2,2]_{\mathfrak a}
\end{align*}
plus all permutations.
\end{itemize}
\end{thm}

The result is proven using the procedure given in \cite{Bi2}.  Since this is
straightforward, we will not give the details.  
It turns out that there are no inner walls for this case either, though this is
certainly tedious to check by hand.
For example, take
$[1,2,3]_{\mathfrak a}=1$ for the $I=3$ case. 
This is compatible with the first set of inequalities.  However, since the total
index is $3$, we have
$[3,3,2]_{\mathfrak a}+[2,1,1]_{\mathfrak a}=2$, and this violates the inequality
$[3,3,2]_{\mathfrak a}\geq 2$.

Indeed, by combining Propositions \ref{P:structure} (3) and \ref{maximalrank1}, 
  and using the connectivity of the 
moduli of parabolic bundles, one can show that the smallest
$U(n)$ for which inner walls can appear is $n=5$ (still assuming $\ell=3$).

\medskip
\noindent \emph{Acknowledgments.}
Proposition \ref{P:isotropic} below has been independently proven by
Florent Schaffhauser in \cite{S}.  His method is to realize  the Lagrangian
representations as fixed points of an antisymplectic involution acting on the space of all unitary representations.
The authors would like to thank him
for many discussions about this problem.
They are also grateful to the mathematics departments at Johns Hopkins University 
 and the  Universit\'e Paris VI for their generous hospitality during the
course of this research.  Funding for this work was provided by 
a US/France Cooperative Research Grant: NSF OISE-0232724, CNRS 14551.
R.W. received additional support  from NSF DMS-9971860.


\bigskip


\frenchspacing

\begin{thebibliography}{99}


\bibitem[AMM]{AMM} A. Alekseev, A. Malkin and E. Meinrenken, 
Lie group valued moment maps.  J. Diff. Geom. 48(3) (1998), 445--495.



\bibitem[AW]{AW} S. Agnihotri and C. Woodward, 
Eigenvalues of products of unitary matrices and quantum Schubert calculus. 
 Math. Res. Lett. 5 (1998), no. 6, 817--836.

\bibitem[Be]{Be} P. Belkale, 
Local systems on $\mathbb P\sp 1-S$ for $S$ a finite set. Compositio Math. 129 
(2001), no. 1, 67--86. 

\bibitem[Bi1]{Bi1} I. Biswas,  A criterion for the existence of a parabolic stable
bundle of rank two over the projective line.
 Int. J. Math. 9 (1998), no. 5, 523--533.

\bibitem[Bi2]{Bi2}  I. Biswas, On the existence of unitary flat connections over 
the punctured sphere with given local monodromy around 
the punctures. Asian J. Math. 3 (1999), no. 2, 333--344.

\bibitem[FMS]{FMS} E. Falbel, J.-P. Marco, and 
F. Schaffhauser, Classifying triples of Lagrangians
in a hermitian vector space, 
 Topology  Appl. 144 (2004), 1--27.

\bibitem[F]{F} W. Fulton, 
Eigenvalues, invariant factors, highest weights, and Schubert Calculus.
Bull. Amer. Math. Soc. 37 (2000), no. 3,  209--249.

\bibitem[G1]{G1} W. Goldman, The symplectic nature of fundamental groups of surfaces, Adv. Math. 54 (1984), 200--225.

\bibitem[G2]{G2} W. Goldman, Invariant functions on Lie groups and Hamiltonian flows of surface group representations, Invent. Math. 85 (1986), 1--40.

\bibitem[GM]{GM} W. Goldman and J. Millson, Eichler-Shimura homology and the finite generation of cusp forms by hyperbolic Poincar\'e series, Duke Math. J. 53 (1986), 1081--1091.

\bibitem[JW]{JW} L. Jeffrey and J. Weitsman, Bohr-Sommerfeld orbits in the moduli
space of flat connections and the Verlinde dimension formula. Comm. Math. Phys. 150
(1992), 593--630.

\bibitem[JM]{JM} D. Johnson and J. Millson, Deformation spaces associated to compact hyperbolic manifolds.  Discrete groups in geometry and analysis (New Haven, Conn., 1984),  48--106, Progr. Math., 67, BirkhŠuser Boston, Boston, MA, 1987. 

\bibitem[K]{K} A. Klyachko, Stable bundles, 
representation theory and Hermitian operators. 
Selecta Math. (N.S.) 4 (1998), no. 3, 419--445.

\bibitem[KM]{KM} M. Kapovich and J. Millson, 
The symplectic geometry of polygons in Euclidean space.
J. Differential Geom. 44 (1996), no. 3, 479--513.

\bibitem[KS]{KS} M. Kashiwara and P. Schapira, 
``Sheaves on Manifolds," Grundlehren der Mathematischen Wissenshaften 292,
Springer-Verlag, Berlin, 1990.

\bibitem[LM]{LM} A. Lubotzky and A. Magid, Varieties of representations of finitely generated groups, Mem. Amer. Math. Soc. 58 (1985), no. 336.

\bibitem[MS]{MS} V. Mehta and C. Seshadri, Moduli of vector bundles on curves
with parabolic structures. Math. Ann. 248 (1980), no. 3, 205--239. 

\bibitem[NS]{NS} M. Narasimhan and C. Seshadri, 
Stable and unitary vector bundles on a compact Riemann surface.
Ann. Math. (2) 82  (1965), 540--567. 

\bibitem[S]{S} F. Schaffhauser, Representations of the fundamental group
of an $l$-punctured sphere generated by products of Lagrangian involutions,
to appear in Can. Math. J.




\bibitem[W]{W} A. Weil,
Remarks on the cohomology of groups.
Ann. Math. (2) 80  (1964), 149--157. 


\end{thebibliography}
\end{document}